\documentclass[12pt,reqno]{amsart}
\usepackage[textwidth=468pt,textheight=648pt,centering]{geometry}
\usepackage[T1]{fontenc}
\usepackage[utf8]{inputenc}
\usepackage{amsmath,amssymb,amsthm,amscd,mathtools}
\usepackage{graphicx,microtype,enumitem,float,needspace}
\usepackage{tikz}
\usetikzlibrary{arrows.meta}
\usepackage[hidelinks]{hyperref}
\hypersetup{pdftitle={Meridional ranks of Whitehead doubles, Part I}}
\usepackage{etoolbox,xparse}

\usepackage[
  backend=biber,
  style=numeric,
  sorting=nyt,
  giveninits=true,
  doi=false,
  url=false,
  isbn=false,
  maxnames=10
]{biblatex}

\usepackage{adjustbox}

\ExplSyntaxOn
\NewExpandableDocumentCommand{\Z}{m}
  {\prg_replicate:nn {#1} {0 &}}
\ExplSyntaxOff

\DeclareFieldFormat[article,incollection,inproceedings,misc,unpublished]{title}{\mkbibemph{#1}}
\DeclareFieldFormat{journaltitle}{#1}
\DeclareFieldFormat{journalsubtitle}{#1}
\DeclareFieldFormat[article,incollection,inproceedings,misc]{volume}{\mkbibbold{#1}}
\DeclareFieldFormat{pages}{#1}
\DeclareFieldFormat[article]{number}{\ifstrequal{#1}{Special Issue}{special issue}{no.\addnbspace#1}}
\renewbibmacro*{journal+issuetitle}{%
  \usebibmacro{journal}%
  \setunit*{\addspace}%
  \iffieldundef{series}{}{\printfield{series}\setunit{\addspace}}%
  \printfield{volume}%
  \setunit{\addspace}%
  \usebibmacro{issue+date}%
  \setunit{\addcomma\space}%
  \printfield{number}%
  \setunit{\addcomma\space}%
  \printfield{eid}%
  \newunit}
\DeclareFieldFormat{eprint:arxiv}{%
  arXiv\addcolon\space\mbox{\ifhyperref
    {\href{https://arxiv.org/abs/#1}{\nolinkurl{#1}}}
    {\nolinkurl{#1}}}}
\renewbibmacro{in:}{}
\AtEveryBibitem{\ifentrytype{book}{\clearfield{pages}}{}}
\AtEveryBibitem{\clearfield{eprintclass}}
\AtEveryBibitem{\clearlist{language}}
\AtEveryBibitem{\ifentrytype{article}{\clearfield{eprint}\clearfield{eprinttype}\clearfield{month}\clearfield{day}}{}}

\let\OriginalSubsection\subsection
\RenewDocumentCommand{\subsection}{s o m}{%
  \IfBooleanTF{#1}{\OriginalSubsection*{#3}}{%
    \IfNoValueTF{#2}{%
      \ifstrequal{#3}{$SL_2(\mathbb{C})$ character variety of knots}
        {\OriginalSubsection{\texorpdfstring{#3}{SL2(C) character variety of knots}}}
        {\ifstrequal{#3}{$PSL_2(\mathbb{C})$ representation variety}
          {\OriginalSubsection{\texorpdfstring{#3}{PSL2(C) representation variety}}}
          {\OriginalSubsection{#3}}}%
    }{\OriginalSubsection[#2]{#3}}%
  }%
}
\graphicspath{{figures/}}
\theoremstyle{plain}
\newtheorem{theorem}{Theorem}[section]
\newtheorem{proposition}[theorem]{Proposition}
\newtheorem{lemma}[theorem]{Lemma}
\newtheorem{corollary}[theorem]{Corollary}

\theoremstyle{definition}
\newtheorem{definition}[theorem]{Definition}
\newtheorem{example}[theorem]{Example}

\newtheorem{remark}[theorem]{Remark}
\numberwithin{equation}{section}
\newcommand{\C}{\mathbb C}

\newcommand{\Wh}{\operatorname{Wh}}

\newcommand{\Span}{\operatorname{Span}}
\newcommand{\ar}{\operatorname{ar}}
\newcommand{\mr}{\operatorname{mr}}
\newcommand{\gr}{\operatorname{gr}}
\newcommand{\im}{\operatorname{im}}

\newcommand{\Cord}{\operatorname{Cord}}
\newcommand{\GL}{\operatorname{GL}}
\newcommand{\SL}{\operatorname{SL}}
\newcommand{\PSL}{\operatorname{PSL}}
\newcommand{\SO}{\operatorname{SO}}
\newcommand{\murk}{\mu\operatorname{rk}}

\newcommand{\diag}{\operatorname{diag}}
\newcommand{\bdefn}{\begin{definition}}
\newcommand{\edefn}{\end{definition}}
\newcommand{\blem}{\begin{lemma}}
\newcommand{\elem}{\end{lemma}}
\newcommand{\bpf}{\begin{proof}}
\newcommand{\epf}{\end{proof}}
\newcommand{\brem}{\begin{remark}}
\newcommand{\erem}{\end{remark}}
\newcommand{\aand}{\ \text{and}\ }

\definecolor{whtubeedge}{RGB}{99,145,166}
\definecolor{whtubefill}{RGB}{225,239,244}
\definecolor{whcoreblue}{RGB}{20,78,105}
\title[Meridional ranks of iterated Whitehead doubles]{Meridional ranks of iterated Whitehead doubles}
\author[H. Gao]{Honghao Gao}
\address[H. Gao]{Yau Mathematical Sciences Center, Tsinghua University, Beijing, 100084 \& Beijing Institute of Mathematical
Sciences and Applications, Beijing, 101408, China}
\email{gaohonghao@tsinghua.edu.cn}
\author[S. Gu]{Shijie Gu}
\address[S. Gu]{Department of Mathematics, Northeastern University, Shenyang, Liaoning, China, 110004}
\email{shijiegutop@gmail.com}
\author[Z. Li]{Zhenkun Li}
\address[Z. Li]{State Key Laboratory of Mathematical Sciences,
Academy of Mathematics and Systems Science, Chinese Academy of Sciences,
Beijing 100190, China}
\email{zhenkun@amss.ac.cn}
\author[J. Wang]{Jian Wang}
\address[J. Wang]{State Key Laboratory of Mathematical Sciences,
Academy of Mathematics and Systems Science, Chinese Academy of Sciences,
Beijing 100190, China}
\email{jian.wang.4@amss.ac.cn}
\date{\today}
\begin{document}
\begin{abstract}
We prove that every iterated untwisted Whitehead double of a nontrivial two-bridge knot satisfies the meridional rank conjecture. More precisely, for such a knot $K$ and every $N\geq0$, the augmentation rank, meridional rank, and bridge number of $\Wh^N(K)$ are all $2^{N+1}$. We use microlocal sheaf theory and a doubling construction to produce a lower bound for the meridional rank of the Whitehead double, assuming the peripheral data satisfy suitable conditions. For two-bridge knots, the required curves arise from ideal points detecting nonzero boundary slopes. In contrast, the ordinary group rank of an $N$-fold Whitehead double is at most the rank of the companion group plus $N$. We also give a non-example when the doubling construction does not work -- the first Whitehead double of the square knot has augmentation rank at most 5, whereas its meridional rank and bridge number are 6.
\end{abstract}
\maketitle
\tableofcontents
%\clearpage
% BEGIN INCLUDED SOURCE: sections/01_introduction.tex
\section{Introduction}\label{sec:introduction}

Let $K\subset S^3$ be a knot and let $\pi_K=\pi_1(S^3\setminus K)$. The \emph{meridional rank} $\mr(K)$ is the minimum number of meridians needed to generate $\pi_K$. A bridge presentation gives a meridional generating set, so $\mr(K)\leq b(K)$, where $b(K)$ denotes the bridge number. The meridional rank conjecture of Cappell--Shaneson asks whether equality holds for every knot. It appears in Kirby's problem list \cite[Problem 1.11]{Kirby1997} and, in its link formulation. See also \cite[Problem 1.18(a)]{K3}.

The conjecture is known for several classes of knots and links. Boileau--Zieschang proved it for Montesinos links by using their twofold branched covers \cite{BoileauZieschang}, and Rost--Zieschang proved it for torus links \cite{RostZieschang}. The characterization of knots whose groups can be generated by two meridians is due to Boileau--Zimmermann \cite[Corollary 3.3]{BoileauZimmermann}, together with the subsequent orbifold geometrization theorem \cite{BoileauLeebPorti}.

Diagrammatic descriptions of bridge number clarify a difficulty in the conjecture. Blair--Kjuchukova--Velazquez--Villanueva \cite{BlairKjuchukovaVelazquezVillanueva} proved that the bridge number equals the Wirtinger number, the minimum number of seed strands needed for successive Wirtinger generation over all diagrams. Blair--Kjuchukova--Pfaff allow relations witnessed by certain embedded two-spheres; their plain sphere number also equals the bridge number \cite{BlairKjuchukovaPfaffPlain}. To deduce the conjecture from these descriptions, one would still need to pass from an arbitrary meridional generating set to such a system without increasing its size.

For lower bounds, the conjugacy class of a meridian is essential. Ordinary group rank $\gr(K)$ allows arbitrary generators and can be much smaller: a nontrivial torus-knot group has rank two, whereas its meridional rank can be arbitrarily large \cite{RostZieschang}. A quotient retaining the meridian class can give a stronger bound. Baader--Blair--Kjuchukova use meridian-to-reflection Coxeter quotients to prove the conjecture for twisted links and arborescent links associated to bipartite trees with nonzero even weights \cite{BaaderBlairKjuchukova}. Baader--Kjuchukova \cite[Section~2]{BaaderKjuchukova} use symmetric-group quotients to prove the meridional rank conjecture for a family of twisted Whitehead doubles of connected sums of $T(2,5)$, with arbitrarily large bridge number. Recent work of Blair--Kjuchukova--Pfaff extends this approach to fishnet families with any prescribed bridge number at least three and arbitrarily large bridge distance relative to that number \cite{BlairKjuchukovaPfaffFishnet}. The question of what finite quotients can detect is also recorded in \cite[Problem 1.18(b)]{K3}.

Knot contact homology gives a different lower bound. Cornwell \cite[Theorems~1.2 and~1.5]{Cornwell} relates
its augmentations to representations of the knot group and proves
\begin{equation}\label{eq:intro-rank-chain}
 \ar(K)\leq\mr(K)\leq b(K),
\end{equation}
where $\ar(K)$ denotes augmentation rank. In the corresponding microlocal rank one representations, a meridian acts by a diagonalizable rank-one perturbation of the identity. For an irreducible microlocal rank one
representation $\rho$, the images of $\rho(m_i)-I$ for any meridional generating set $\{m_i\}$ span the representation space. Its dimension, therefore, bounds the meridional rank from below. We will vary these representations and their peripheral eigenvalues along an algebraic curve to choose data for successive Whitehead doubles. We recall the definitions in Section~\ref{subsec:augmentation-representations}.

We study these ranks under untwisted Whitehead doubling. Two satellite results are particularly relevant. For satellites with braid patterns, Cornwell--Hemminger prove a maximal augmentation-rank multiplication theorem and apply it to the iterated torus knots satisfying their hypotheses \cite[Theorem 1.2 and Corollary 1.3]{CornwellHemminger}. The Whitehead pattern has winding number zero, so that theorem does not apply directly. Using the action on the tree associated to a satellite decomposition, Boileau--Dutra--Jang--Weidmann control the peripheral intersections of meridional subgroups. They prove the conjecture for the class generated from meridionally tame knots by connected sums and satellites with braid patterns \cite[Theorem 1]{BDJW}. This includes every knot whose exterior is a graph manifold \cite[Corollary 1]{BDJW}. However, Whitehead doubles of nontrivial knots are not meridionally tame \cite[Introduction and Section 6]{BDJW}, so this theorem does not directly include Whitehead iteration.

Dutra's preprint \cite[Theorem 1.1]{Dutra} states the meridional rank conjecture for (once) Whitehead doubles of prime algebraically tame knots. Even with this conclusion, iteration would require a new argument: the first Whitehead double is no longer algebraically tame \cite[Example 2.8]{Dutra}.\footnote{There is also an issue in the
pattern-subgroup argument of arXiv:2212.13081v1. In the notation on p.~7 of that version, $m_V=x_2x_1^{-1}x_2^{-1}x_1=\alpha(y_1y_2)=\omega(y_1y_2)$ is a nontrivial element of $U_\alpha\cap U_\omega$, contradicting Lemma 3.1 with $g=1$. That lemma is used in the proof of Lemma 3.8(2). This identifies a gap in the argument, not a counterexample to the conclusion of Theorem 1.1. Our proof does not use these subgroup lemmas.} Our contribution is a sheaf-level doubling construction whose iteration is controlled on peripheral curves realized by the companion. By keeping track of both peripheral eigenvalues, we can choose representations for which every required doubling step is defined.

Throughout, $\Wh(K)$ denotes the untwisted Whitehead double with the pattern and peripheral conventions fixed in Section~\ref{subsec:whitehead-pattern}. We set $\Wh^0(K)=K$ and $\Wh^{N+1}(K)=\Wh(\Wh^N(K))$. For a nontrivial companion,
\begin{equation}\label{eq:intro-bridge-doubling}
 b(\Wh(K))=2b(K);
\end{equation}
see \cite[Section 4]{Dutra}. Thus, an augmentation rank lower bound that doubles at every step proves the meridional rank conjecture for the iterates, provided the initial augmentation rank equals $b(K)$.

\Needspace{9\baselineskip}
\subsection{Main results}

Our first main result applies to every nontrivial two-bridge knot.

\begin{theorem}\label{thm:main-two-bridge}
Let $K$ be a nontrivial two-bridge knot. For every integer $N\geq0$,
\[
 \ar(\Wh^N(K))=\mr(\Wh^N(K))=b(\Wh^N(K))=2^{N+1}.
\]
\end{theorem}

For $N=1$, the meridional-rank equality already follows from the lower bound for prime satellite knots in \cite[Corollary 1.6]{BJW}; see also \cite[Introduction]{BDJW}. Theorem~\ref{thm:main-two-bridge} also provides maximal-rank augmentations and propagates both equalities through every finite number of Whitehead doubles.

Given an irreducible microlocal rank one representation $\rho$ on $V$, we extend $\rho\oplus\rho^*$ from the companion group to the satellite group, where $(V^*,\rho^*)$ is the dual representation. The two distinguished peripheral eigenlines span a two-dimensional subspace of $V\oplus V^*$, on which we solve the Whitehead-pattern relations explicitly. The remaining subspace is fixed by the two pattern meridians. Under the nonvanishing conditions of Theorem~\ref{thm:one-step-doubling}, this gives an irreducible $2\dim V$-dimensional representation whose meridian still has microlocal rank $1$, doubling the resulting lower bound for meridional rank.

To iterate, we compute the output longitude eigenvalue as well, so that the output of one step can be tested against the hypotheses of the next. This gives the algebraic correspondence of Section~\ref{subsec:iteration}. We must then find peripheral values realized by a fixed companion that admits any prescribed finite number of steps.

The difficulty is that the peripheral values realized by a fixed knot might all lie in the exceptional locus. An open dense set in the full peripheral torus does not rule this out. We address this issue using a branch at infinity of the curve. The branch is called \emph{Whitehead-iterable branch} if the meridian eigenvalue has a pole and the squared longitude eigenvalue either has a pole or has a zero of order greater than that pole order; see Definition~\ref{def:escaping-branch}. The key observation is that one branch of the doubling correspondence again satisfies this condition. By Theorem~\ref{thm:escaping-nontrapping}, for each $N\geq0$ there is a nonempty Zariski-open subset of the curve on which $N$ doubling steps can be performed. For two-bridge knots, the ideal-point argument in Section~\ref{sec:slope-detection} supplies a curve with an iterable branch whose peripheral values are realized by rank-two augmentations on a Zariski-dense subset. Combining these facts with \eqref{eq:intro-rank-chain} and \eqref{eq:intro-bridge-doubling} proves Theorem~\ref{thm:main-two-bridge}.

The same construction provides a comparison with ordinary group rank. We call an irreducible curve in the peripheral torus an \emph{actual rank-$d$ augmentation peripheral curve} if a Zariski-dense subset of its points is realized by rank-$d$ augmentations with $\epsilon(\mu)\neq1$. The coordinates are $(\mu,\lambda^2)$, where $\lambda$ is the preferred longitude eigenvalue; see Definition~\ref{def:actual-kch-curve}.

\begin{theorem}\label{thm:main-ordinary}
For every knot $K$ and every integer $N\geq0$,
\[
 \gr(\Wh^N(K))\leq\gr(K)+N.
\]
If $K$ has an actual rank-$d$ augmentation peripheral curve with an iterable branch, then
\[
 \ar(\Wh^N(K))\geq2^Nd
\]
for every $N\geq0$. In particular,
\[
 \gr(\Wh^N(K))\leq\ar(\Wh^N(K))
 \quad\text{whenever}\quad \gr(K)+N\leq2^Nd.
\]
For a nontrivial two-bridge knot, this gives
\[
 \gr(\Wh^N(K))\leq N+2\leq
 \ar(\Wh^N(K))=2^{N+1}.
\]
\end{theorem}

The proof is given in Section~\ref{sec:ordinary-rank}. A complementary lower bound for ordinary group rank was obtained by Gu--Wang--Zou
\cite[Proposition 1.3]{GuWangZou}. For every nontrivial knot $K$ and $N\geq1$, their estimate and Theorem~\ref{thm:main-ordinary} give
\[
 N+1\leq\gr(\Wh^N(K))\leq N+\gr(K).
\]
Thus ordinary group rank grows linearly, whereas the augmentation and meridional ranks in Theorem~\ref{thm:main-two-bridge} grow exponentially. The lower bound in \cite{GuWangZou} uses the hyperbolic pieces of the JSJ decomposition and holds without an augmentation hypothesis. That work also determines ordinary group rank exactly for hyperbolic two-bridge companions; see Section~\ref{sec:ordinary-rank}.

The augmentation lower bound is not always sharp; this limitation is already present in Cornwell's work \cite{Cornwell}. Our last main result exhibits a further obstruction to iteration: even when the augmentation rank agrees with the bridge number for the companion, it may fail to double. The example nevertheless satisfies the meridional rank conjecture.

\begin{theorem}\label{thm:main-square}
Let $S=T(2,3)\#\overline{T(2,3)}$ be the square knot. Then
\[
 \ar(S)=\mr(S)=b(S)=3,
\]
whereas
\[
 \ar(\Wh(S))\leq5<6=\mr(\Wh(S))=b(\Wh(S)).
\]
\end{theorem}

The augmentation obstruction uses the two trefoil factors of the companion group and their preferred longitudes. The meridional lower
bound uses a different representation: an explicit epimorphism $\pi_{\Wh(S)}\to S_{14}$ and its permutation module $\mathbb F_2^{14}$.
The first cohomology of $\pi_{\Wh(S)}$ with coefficients in this module has dimension $14$, while the meridian defect has rank $3$. A cocycle
estimate then rules out generation by five meridians; see Section~\ref{sec:square-meridional} for the representation and exact
linear calculation. Here the finite quotient supplies the coefficient module for the cohomology of the knot group. This gives a strict gap
between augmentation rank and meridional rank. The exact augmentation rank of $\Wh(S)$ and the ranks of further Whitehead doubles of $S$
remain undetermined in this paper.

\subsection{Conventions and organization}

All knots are oriented. A meridian is positively oriented, and elements represented by meridians are understood up to conjugacy. We use a preferred zero-linking longitude. Sheaves and augmentations may be defined over any field $\Bbbk$, but in the algebraic-curve applications, we take $\Bbbk=\mathbb C$. The finite-dimensional permutation module in Section~\ref{sec:square-meridional} is explicitly over $\mathbb F_2$.

We use $\pi_K :=\pi_1(S^3\setminus K)$ for the knot group throughout. The commutator convention is $[a,b]=aba^{-1}b^{-1}$. Products of permutations are composed from right to left, and permutation matrices act on column vectors by $g e_i=e_{g(i)}$. The embedding $\phi_K$ of the Whitehead solid torus into a tubular neighborhood of $K$ respects the preferred zero framing. A ``two-bridge knot'' in the statements below is nontrivial.

Section~2 recalls the representations and augmentations and constructs the doubling correspondence. Section~\ref{sec:meridional-bounds} gives
the meridional upper bounds. Section~\ref{sec:algebraic-preliminaries} recalls the algebraic curves, valuations, and character varieties used to find peripheral data to which the iteration applies. Section~\ref{sec:escaping-curves} proves the iterable-branch criterion, and Section~\ref{sec:slope-detection} applies slope detection to prove Theorem~\ref{thm:main-two-bridge}. Section~\ref{sec:ordinary-rank}
compares ordinary group rank with augmentation rank. Finally, Section~\ref{sec:first-doubles} proves Theorem~\ref{thm:main-square} for
the first Whitehead double of the square knot.
% END INCLUDED SOURCE: sections/01_introduction.tex
% BEGIN INCLUDED SOURCE: sections/02_doubling.tex
\section{Microlocal rank-1 representations and Whitehead doubling}\label{sec:shead_aug}

The meridian and the preferred longitude of a knot $K$ are denoted by $m$ and $\ell$, respectively.  We use the specialization of knot contact homology in which the extra variable $U$ is set equal to one.  For an embedded submanifold $A\subset S^3$, we denote by $\mathcal{O}p(A)$ a fixed sufficiently small open tubular neighborhood of $A$.

\subsection{Microlocal rank-1 representations and augmentations}
\label{subsec:augmentation-representations}

We first introduce the notation for the microlocal rank of a representation of the knot group.

\begin{definition}
\label{def:microlocal-rank}
For a representation $\rho\colon\pi_K\to\GL(V)$, define the microlocal rank to be
\[
 \murk_K(\rho):=\operatorname{rank}\bigl(\rho(m)-I\bigr),
\]
where $I$ is the identity matrix of size $\dim V$. This number is independent of the chosen meridian, since all meridians are conjugate. 

We call $\rho$ a \emph{$\murk$-one representation} if
$\murk_K(\rho)=1$.

A $\murk$-one representation is \emph{diagonalizable} if there is a basis in which
\[
\rho(m)\sim\diag(\mu,1,\ldots,1),\qquad \mu\neq1.
\]
A $\murk$-one representation is \emph{unipotent} if there is a basis in which
\[
 \rho(m)\sim
 \begin{pmatrix} 1& 1 \\ 0& 1\end{pmatrix}\oplus I_{\dim V -2}.
\]
\end{definition}
 
We recall the sheaf-theoretic reason for this terminology \cite[Section~3.2]{GaoSimple}.  Let $j\colon S^3\setminus K\hookrightarrow S^3$ be the open embedding. A representation $\rho\colon\pi_K\to\GL(V)$ determines a local system $\mathcal E_\rho$ on the knot complement. For the underived direct image, the stalk of $j_*\mathcal E_\rho$ on $K$ is the meridian-invariant subspace $V^{\rho(m)}$. The local gluing map is the inclusion
\[
 T\colon V^{\rho(m)}\hookrightarrow V.
\]
Its mapping cone is the microlocal stalk, up to a degree shift, in the sense of Kashiwara--Schapira, along the knot conormal torus $\Lambda_K$ inside the contact cosphere bundle, i.e. $\Lambda_K\subset ST^*S^3$ is a Legendrian surface.  Consequently,
\[
 \dim H^\bullet\!\left(\operatorname{Cone}(T)\right)
   =\dim(V/V^{\rho(m)})
   =\operatorname{rank}\bigl(\rho(m)-I\bigr).  
\]
Thus a $\murk$-one representation induces a microlocal rank one sheaf microsupported along the conormal. Conversely, up to the equivalence modulo local systems on $S^3$ in \cite[Definition~3.8 and Proposition~3.11]{GaoSimple}, a microlocal rank one abelian sheaf comes from one of the three forms: (1) an irreducible diagonalizable $\murk$-1 representation via the above construction; (2) an irreducible unipotent $\murk$-1 representation via the above construction; and (3) a local system of rank $1$ that is supported on $K$ \cite{GaoSimple}. The first type is generic and the main focus of this paper.

A diagonalizable $\murk$-one representation is also called a KCH representation by Cornwell \cite{Cornwell}, where KCH refers to \emph{knot contact homology}. Next, we explain its background related to knot contact homology and augmentations.

The (framed) knot contact homology is a knot invariant defined via counting holomorphic disks for $\Lambda_K$, whose formula was first discovered by Ng \cite{Ng}. We use the specialization $U=1$ of the Legendrian dg algebra $\mathcal A(\Lambda_K)$ associated to $\Lambda_K$. The dg algebra $\mathcal A(\Lambda_K)$ is semi-freely generated by Reeb chords over $\mathbb Z[H_1(\Lambda_K)]\cong \mathbb Z[\mu^{\pm1},\lambda^{\pm1}]$. An augmentation of the dg algebra is a morphism $\epsilon: \mathcal A(\Lambda_K)\rightarrow \Bbbk$ which preserves grading, differential, and the unit. The case of the knot conormal has the special property that its generators are concentrated in non-negative degrees. Hence an augmentation descends to the degree $0$ homology $\epsilon:  H_0(\mathcal A(\Lambda_K))\rightarrow \Bbbk$.

The algebra $H_0(\mathcal A(\Lambda_K))$ admits a combinatorial model by the framed cord algebra.
\begin{definition}
\label{def:cord-algebra}
The \emph{framed cord algebra} $\Cord(K)$ is the tensor algebra over $\mathbb Z[\lambda^{\pm1},\mu^{\pm1}] $ generated by symbols $[\gamma]$, $\gamma\in\pi_K$, modulo the following relations:
\begin{enumerate}
	\item $[e]=1-\mu$,                               
	\item $[m\gamma]=[\gamma m]=\mu[\gamma]$,  $[\ell\gamma]=[\gamma\ell]=\lambda[\gamma]$,
	\item  $[\gamma_1\gamma_2]-[\gamma_1m\gamma_2] =[\gamma_1][\gamma_2]$.  
\end{enumerate}

An \emph{augmentation} is a unital algebra homomorphism
\[
 \epsilon\colon\Cord(K)\longrightarrow\Bbbk.
\]
Note that $\epsilon$ includes a specialization of $\lambda$ and $\mu$ to nonzero elements of $\Bbbk$.
\end{definition}

Augmentations of the framed cord algebra and the above equivalence classes of microlocal rank one abelian sheaves modulo local systems are in bijective correspondence, with the normalization of \cite[Theorem~4.17]{GaoSimple}. This is a theorem under the principle of ``augmentation-sheaf correspondence'' in the theory of Legendrian invariants. Considering the first and second cases of the sheaf classification, suppose $(\rho, V)$ is a $\murk$-1 representation; then the associated augmentation $\epsilon_\rho$ is given by
\[
\epsilon_\rho([\gamma]) = \operatorname{tr}(\rho(\gamma) - \rho(m\gamma)).
\]

Conversely, suppose $\epsilon$ is an augmentation of the framed cord algebra with $\epsilon(\mu)\neq1$. To reconstruct the $\murk$-1 representation $(\rho_\epsilon, V_\epsilon)$ following \cite[Proposition~4.11]{GaoSimple}, we choose a set of based meridional generators $\{m_1,\ldots,m_N\}$ with $m_1=m$. Since the meridians are conjugate, there exist elements $g_i\in\pi_K$ such that $ m_i=g_i^{-1}m_1g_i$. Define an $N$ by $N$ square matrix $R_\epsilon$ by $(R_\epsilon)_{ij}:=\epsilon([g_ig_j^{-1}])$. Let $R_j$, $j=1,\ldots,N$, be the column vectors of $R_\epsilon$, and define the vector space
\begin{equation}
 V_\epsilon:=\Span_{\Bbbk}\{R_1,\ldots,R_N\}.
 \label{eq:augmentation-matrix}
\end{equation}
and for any $\gamma\in \pi_K$, the $\rho_\epsilon(\gamma)$ action is given by 
\[
 \rho_\epsilon(\gamma)R_j
   :=\bigl(\epsilon([g_i\gamma g_j^{-1}])\bigr)_{i=1}^N.
\]
The resulting representation is indeed $\murk$-1. For example, by relation 3 in Definition \ref{def:cord-algebra},
\begin{equation*}
 \rho_\epsilon(m_t)R_j
   =R_j-\epsilon([g_tg_j^{-1}])R_t.
\end{equation*}
It follows at once that $\operatorname{rank}(\rho_\epsilon(m_t)-I)\leq1$. Since $(R_\epsilon)_{tt}=1-\epsilon(\mu)\neq0$, we have $R_t\neq0$ and $\rho_\epsilon(m_t)R_t=\epsilon(\mu)R_t$. Thus the rank is exactly one, and the representation is diagonalizable $\murk$-one. This elementary identity is the point at which the microlocal condition enters the construction.

The dimension of $V_\epsilon$ is independent of $N$ and the choices of meridional generators. Meanwhile, it does depend on the choice of $\epsilon$. Moreover, the representation is irreducible \cite[Theorem~1.2]{Cornwell}.

\begin{definition}
For a knot $K$, define its \emph{augmentation rank}, \emph{meridional rank}, and
\emph{bridge number} as follows:
\begin{align*}
 \operatorname{ar}(K)
   &:=\max\{\operatorname{dim}V_\epsilon:
       \epsilon\colon\Cord(K)\to\Bbbk,\ \epsilon(\mu)\neq1\},\\
 \operatorname{mr}(K)
   &:=\min\{r:\pi_K\text{ is generated by }r\text{ meridians}\},\\
 b(K)
   &:=\min\{r:K\text{ admits an }r\text{-bridge presentation}\}.
\end{align*}
\end{definition}

Via the augmentation-sheaf correspondence, the augmentation rank can also be read as the maximum rank of the irreducible diagonalizable microlocal rank one representations of $\pi_K$. The two occurrences of the word rank should not be confused. Microlocal rank records the normal jump along the conormal, whereas the other rank is rank in the usual sense.

We immediately have the following inequalities.
\begin{proposition}
\label{prop:rank-inequalities}
For every knot $K$,
\[
 \operatorname{ar}(K)\leq\operatorname{mr}(K)\leq b(K).
\]
\end{proposition}

\begin{proof}
By construction \eqref{eq:augmentation-matrix}, the maximal value of $\dim V_\epsilon$ is no larger than the number of meridional generators. Hence it yields
\[
 \operatorname{ar}(K)\leq\operatorname{mr}(K).
\]

A bridge presentation with $b(K)$ bridges gives a Wirtinger presentation
whose bridge meridians generate $\pi_K$.  Hence
\[
\operatorname{mr}(K)\leq b(K).
\]
\end{proof}

Let $(\rho,V)$ be a diagonalizable $\murk$-one representation. Since the distinguished meridian and the preferred longitude commute in the knot group, we can choose a suitable basis such that their representing matrices have the form
\[
 \rho(m)=
 \begin{pmatrix}\mu_0&0\\0&I\end{pmatrix},
 \qquad
 \rho(\ell)=
 \begin{pmatrix}\lambda_0&0\\0&*\end{pmatrix},
\]
Note that the distinguished $\mu_0$-eigenvector of $\rho(m)$ is also an eigenvector of $\rho(\ell)$, with eigenvalue $\lambda_0$. This is a property that we shall use repetitively.

Finally, the associated augmentation has 
\begin{equation*}
 \epsilon_\rho(\mu)=\mu_0,\qquad
 \epsilon_\rho(\lambda)=\lambda_0,\qquad
 \epsilon_\rho([\gamma])=(1-\mu_0)\rho(\gamma)_{11}.
\end{equation*}

\begin{definition} Following the notation above, we call the pair
\[
(\mu_0,\lambda_0) \in (\Bbbk^\times)^2
\]
the \emph{peripheral data} of the $\murk$-1 representation $(\rho,V)$ or the associated augmentation.
\end{definition}

\subsection{The Whitehead pattern and the satellite group}
\label{subsec:whitehead-pattern}

Let $\mathrm{Wh}:=C\cup P\subset S^3$ be the oriented Whitehead link in standard position, where $C$ is the axis component and $P$ is the pattern component. Removing a tubular neighborhood of $C$ gives the solid torus
\[
 \mathcal V:=S^3\setminus\mathcal{O}p(C),
\]
and the pattern component becomes a pattern $P\subset\mathcal V$. The pattern has winding number zero. Equivalently, its two longitudinal runs travel once around $\mathcal V$ with opposite orientations. Our peripheral convention on the outer boundary is
\begin{equation*}
 m_C=\ell_{\mathcal V},\qquad \ell_C=m_{\mathcal V},
\end{equation*}
where $(m_C,\ell_C)$ is the peripheral pair of the deleted axis. This fixes the notation in the gluing formulas below.

\begin{definition}
\label{def:whitehead-double}
Let $P\subset\mathcal V$ be the Whitehead pattern fixed above. The \emph{untwisted Whitehead double} of $K$ is the knot
\[
 \Wh(K):=\phi_K(P),
\]
where $\phi_K$ is the zero-framed embedding chosen above. Thus its restriction to the boundary identifies the peripheral curves by
\[
 (\phi_K)_*(m_{\mathcal V})=m_K,\qquad
 (\phi_K)_*(\ell_{\mathcal V})=\ell_K,
\]
with $\ell_K$ the preferred longitude of $K$. The clasp sign is determined by our fixed choice of $P$.
\end{definition}

\begin{figure}[tbhp]
 \centering
 \includegraphics[width=.7\textwidth]{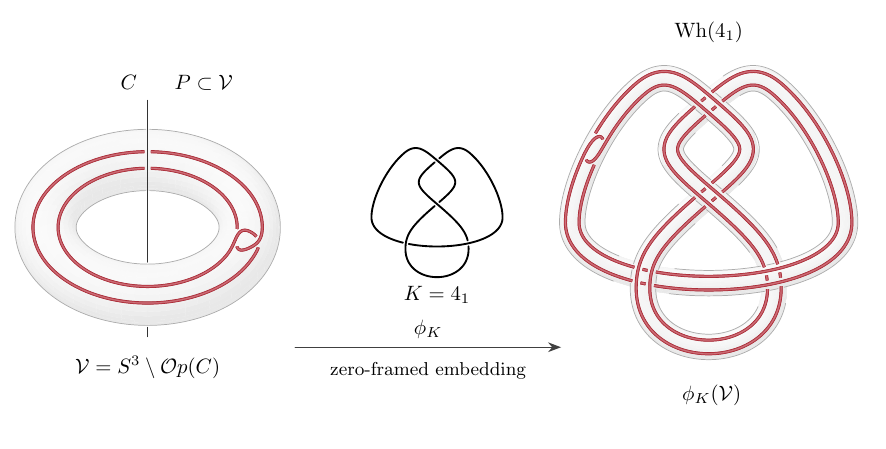}
 \caption{Left: Whitehead pattern. Right: Whitehead double of knot $4_1$.}
 \label{fig:whitehead-fig8-overview}
\end{figure}

Orient all braid strands downward and let $\sigma_i$ denote a positive
right-handed crossing.  One convenient braid representative for the
oriented Whitehead link is
\[
 B_W:=\sigma_1\sigma_2^{-1}\sigma_1\sigma_2^{-1}\sigma_1\in B_3,
\]
where the first strand closes to $C$ and the other two strands close to
$P$.  The four crossings between $C$ and $P$ have signed sum zero.
The unique self-crossing of $P$ accounts for the correction from the
braid framing to the preferred zero framing.

We now pass to the group of the pattern exterior
\[
 E_P:=\mathcal V\setminus\mathcal{O}p(P).
\]
Its fundamental group admits the following presentation \cite{Dutra}:
\begin{equation}
 \pi_1(E_P)=\langle x_1,x_2,\ell_{\mathcal V}\mid\;
 \ell_{\mathcal V}x_2\ell_{\mathcal V}^{-1}=x_1,
 \ell_{\mathcal V}(x_2x_1^{-1}x_2^{-1}x_1x_2^{-1})
 \ell_{\mathcal V}^{-1}=x_2x_1^{-1}x_2^{-1}\rangle .
\label{eq:dutra-presentation}
\end{equation}
The peripheral meridian on the outer boundary is
\begin{equation}
 m_{\mathcal V}=x_2x_1^{-1}x_2^{-1}x_1.               
 \label{eq:outer-meridian}
\end{equation}
and $\langle m_{\mathcal V},\ell_{\mathcal V}\rangle\cong\mathbb Z^2$.
The meridian of the satellite component is $x_1$.

\begin{figure}[tbhp]
 \centering
 \includegraphics[width=.7\textwidth]{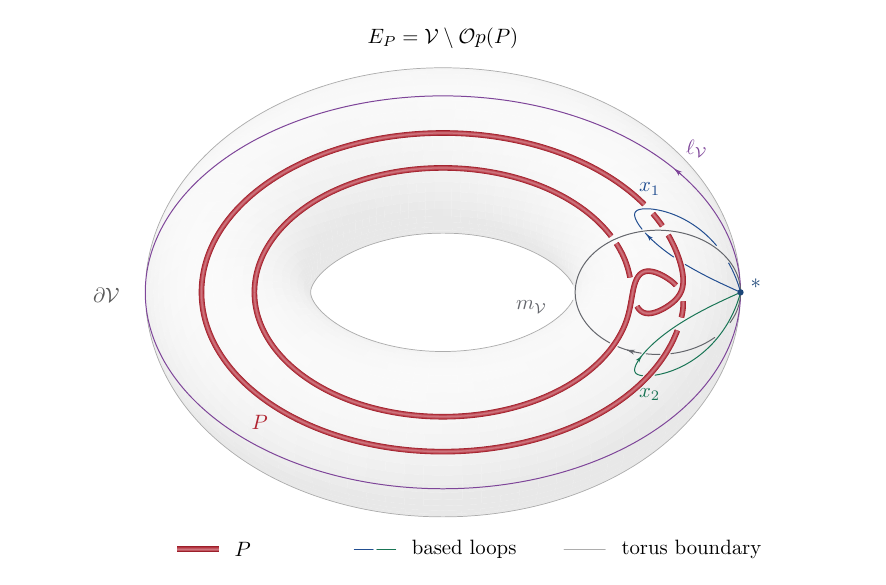}
 \caption{Generating loops of the pattern exterior.}
 \label{fig:whitehead-pattern-generators-closed-longitude}
\end{figure}

\begin{lemma}
\label{lem:pattern-longitude}
With the downward-oriented braid convention fixed above and the
peripheral base path chosen at the meridian $x_1$, the preferred
zero-framed longitude of the satellite component is represented by
\begin{equation}
 \ell_{\Wh}
 =x_1^{-1}\ell_{\mathcal V}x_1\ell_{\mathcal V}^{-1}
  x_1^{-1}\ell_{\mathcal V}^{-1}x_1\ell_{\mathcal V}.
\label{eq:whitehead-longitude}
\end{equation}
\end{lemma}

\begin{proof}
Follow the zero-framed push-off from the basepoint at $x_1$ through the
two oppositely oriented longitudinal runs of the pattern.  The two runs
give the successive conjugations by $\ell_{\mathcal V}$ and
$\ell_{\mathcal V}^{-1}$ in \eqref{eq:whitehead-longitude}.  The four
axis--pattern crossings have total signed contribution zero.  Removing
the single pattern-meridian contribution from the braid framing
therefore gives the displayed preferred longitude, with no additional
power of $m_{\mathcal V}$.
\end{proof}

On the pattern exterior, the torus boundary is subject to the following conditions:
\[
\ell_{\mathcal V}x_2\ell_{\mathcal V}^{-1}=x_1, \qquad
m_{\mathcal V}=x_2x_1^{-1}x_2^{-1}x_1,\qquad
 [\ell_{\mathcal V},m_{\mathcal V}]=1
 \]

\begin{proposition}
\label{prop:van-kampen}
Suppose $K$ is a non-trivial knot. The knot group of the Whitehead double can be computed by 
\begin{equation*}
 \pi_{\Wh(K)}
 \cong \pi_1(E_P)
 *_{\langle m_{\mathcal V},\ell_{\mathcal V}\rangle
      =\langle m_K,\ell_K\rangle}
 \pi_K,
\end{equation*}
where $m_{\mathcal V}$ is identified with $m_K$ and
$\ell_{\mathcal V}$ with $\ell_K$.  
\end{proposition}

\begin{proof}
Let $X_K:=S^3\setminus\mathcal{O}p(K)$. The satellite exterior is the union of $E_P$ and $X_K$.  Their
intersection retracts onto the outer boundary torus of $\mathcal V$.
The zero-framed gluing identifies the two ordered peripheral pairs.
The assertion follows from Seifert–van Kampen.
\end{proof}

\subsection{The one-step doubling construction}
\label{subsec:one-step-doubling}

A $\murk$-one representation of $\pi_{\Wh(K)}$ restricts to a
representation of $\pi_K$ of microlocal rank at most two. The
goal of this subsection is to describe the reverse construction:
starting with a $\murk$-one representation of $\pi_K$, we first pass to
a representation of $\pi_K$ of microlocal rank two and then extend it
across the pattern exterior to obtain a $\murk$-one representation of
$\pi_{\Wh(K)}$.

We first record the algebraic construction used in the doubling
argument.  Let $G$ be a group and let
$\rho\colon G\to\GL(V)$ be a finite-dimensional representation.  Put
$V^*:=\operatorname{Hom}_{\Bbbk}(V,\Bbbk)$.  The dual representation
$(\rho^*,V^*)$ is defined by $
 \bigl(\rho^*(g)\varphi\bigr)(v) :=\varphi\bigl(\rho(g^{-1})v\bigr),
 g\in G,\ \varphi\in V^*,\ v\in V.$

\begin{definition}
\label{def:self-dual-double}
Define the \emph{self-dual double} associated to $\rho$ to be
\begin{equation*}
 \rho^\#:=\rho\oplus\rho^*
 \colon G\longrightarrow\GL(V\oplus V^*).
\end{equation*}
\end{definition}

Let $\rho\colon\pi_K\longrightarrow\GL(V)$ be an $n$-dimensional diagonalizable $\murk$-one representation.  Choose
a distinguished meridian eigenvector $e\in V$.  Since $m_K$ and
$\ell_K$ commute, it is also an eigenvector for the longitude. Let $(\mu,\lambda), \mu\neq1 $ be the peripheral data, namely
\begin{equation*}
 \rho(m_K)e=\mu e, \qquad
 \rho(\ell_K)e=\lambda e.
\end{equation*}

Define $V_0:=\ker(\rho(m_K)-I)$, so that
$V=\Bbbk e\oplus V_0$, and let $e^*\in V^*$ be determined by
$e^*(e):=1$ and $e^*|_{V_0}:=0$.  We consider its self-dual double $(\rho^\#, V\oplus V^*)$ from Definition~\ref{def:self-dual-double}.

Since $\rho(\ell_K)$ preserves $V_0$, the vector $e^*$ is a common
eigenvector for the dual peripheral action, with eigenvalues
$\mu^{-1}$ and $\lambda^{-1}$.  Define 
\begin{equation*}
 U:=\Span\{e,e^*\},\qquad
 W:=\ker\bigl(\rho^\#(m_K)-I\bigr).
\end{equation*}
They form a splitting of the total space, namely:
\begin{equation*}
 V\oplus V^*=U\oplus W.
 \end{equation*}
In the ordered basis $(e,e^*)$ of $U$, put
\begin{equation*}
 M_U:=\rho^\#(m_K)|_U=
 \begin{pmatrix}\mu&0\\0&\mu^{-1}\end{pmatrix},
 \qquad
 L_U:=\rho^\#(\ell_K)|_U= \begin{pmatrix}\lambda&0\\0&\lambda^{-1}\end{pmatrix}.
\end{equation*}
Thus $\rho^\#(m_K)-I$ vanishes on $W$ and has rank two on $U$.
Moreover, both matrices have determinant one. 

We are ready to explain the motivation of the self-dual double construction. 
Under the satellite gluing, the companion meridian $m_K$ is
identified with the commutator $m_{\mathcal V}$ in
\eqref{eq:outer-meridian}.  Consequently, a necessary condition
for a representation of $\pi_K$ to extend over the pattern
exterior is that the image of $m_K$ have determinant one.
The original $\murk$-1 representation does not satisfy this condition,
since $\det\rho(m_K)=\mu\neq1$.  Passing to
$\rho^\#=\rho\oplus\rho^*$ replaces the distinguished eigenvalue
$\mu$ by the reciprocal pair $\mu,\mu^{-1}$, removing this
determinant obstruction.

Our next goal is to extend the self-dual double representation into a $\murk$-1 representation of the Whitehead pattern. Namely, we aim to find $2$-by-$2$ matrices $X_1$ and $X_2$ that satisfy the following boundary conditions:
\begin{equation}
\label{eq:boundary-condition}
L_UX_2L_U^{-1}=X_1,\qquad
 M_U=X_2X_1^{-1}X_2^{-1}X_1.
\end{equation}

Assume $X_1 := \begin{pmatrix}a&b\\ c&d\end{pmatrix}$. The microlocal condition requires $\operatorname{rank}(X_1-I)=1$. Then,
$\det(X_1-I)=0$ gives $bc=(a-1)(d-1)$.  On the chart
$b\neq0$, $d\neq1$, we can normalize the (1,2) entry to $d-1$, and the (2,1) entry
then becomes $a-1$. 

Therefore, we seek a solution of the form
\begin{equation*}
 X_1=
 \begin{pmatrix}
  a&d-1\\
  a-1&d
 \end{pmatrix},
 \qquad
 X_2=L_U^{-1}X_1L_U
 =
 \begin{pmatrix}
  a&\eta^{-1}(d-1)\\
  \eta(a-1)&d
 \end{pmatrix},
\end{equation*}
where $\eta:=\lambda^2$. Let $\widetilde{\mu}:=\det X_1=a+d-1$ be the determinant. For 
$\widetilde{\mu}\notin\{0,1\}$, the matrix $X_1$ is invertible and
diagonalizable $\murk$-1. The same holds for $X_2$, since it is conjugate to $X_1$.
It remains to choose $a,d$ so that the second equation in
\eqref{eq:boundary-condition} holds and $\widetilde{\mu}\notin\{0,1\}$.
For the resulting Whitehead double $\Wh(K)$, $\widetilde{\mu}$ will be
the distinguished meridian eigenvalue.

\begin{lemma}
\label{lem:whitehead-block-existence}
Let $(\mu,\lambda)$ be the peripheral data of the diagonalizable $\murk$-1 representation $\rho$.  Suppose $(\mu,\lambda)$ satisfy the following conditions:
\begin{equation}
 \lambda^2\neq1,\qquad
 \lambda^2\mu\neq-1,\qquad
 \lambda^4\mu^2-\lambda^4\mu+4\lambda^2\mu-\mu+1\neq0.
 \label{eq:doubling-hypotheses}
\end{equation}
Then there exist $a,d\in\Bbbk$ for which
$X_2X_1^{-1}X_2^{-1}X_1=M_U$ and
$\widetilde{\mu}\notin\{0,1\}$.

\end{lemma}

\begin{proof}
We set $r:=(a-1)(d-1)$. Then equation~\eqref{eq:boundary-condition} is equivalent to $X_2X_1^{-1}=M_UX_1^{-1}X_2$. Direct multiplication gives
\begin{align*}
 X_2X_1^{-1}
 &=\frac{1}{\widetilde{\mu}}
 \begin{pmatrix}
  ad-\eta^{-1}r&a(d-1)(\eta^{-1}-1)\\
  d(a-1)(\eta-1)&ad-\eta r
 \end{pmatrix},\\
 M_UX_1^{-1}X_2
 &=
 \frac{1}{\widetilde{\mu}}
 \begin{pmatrix}
  \mu(ad-\eta r)&\mu d(d-1)(\eta^{-1}-1)\\
  \mu^{-1}a(a-1)(\eta-1)&\mu^{-1}(ad-\eta^{-1}r)
 \end{pmatrix}.
\end{align*}
Note that $a\neq1$ and $d\neq1$. Otherwise, say $a=1$ (resp. $d=1$), then $X_1$ and $X_2$ are upper (resp. lower) triangular with the diagonal entries $\{1,d\}$ (resp. $\{a,1\}$), so their commutator has diagonal entries $\{1,1\}$, which is a contradiction to $\mu\neq1$.
Since $\eta\neq1$, comparison of the (1,2) and (1,1) entries
yields
\begin{align}
 a&=\mu d, \label{eq:block-scalar-equations1} \\			
 \eta(1-\mu)ad&+(\mu \eta^2-1)r=0.                                  \label{eq:block-scalar-equations2}
\end{align}

We consider two possibilities depending on the value of $\mu$.

\underline{Case 1, $\mu\neq-1$}.  By $\widetilde{\mu}=a+d-1$ and equation \eqref{eq:block-scalar-equations1}, we have $d=\frac{\widetilde{\mu}+1}{\mu+1}$, and subsequently $
 a=\frac{\mu(\widetilde{\mu}+1)}{\mu+1},
 r=\frac{(\mu\widetilde{\mu}-1)(\widetilde{\mu}-\mu)}{(\mu+1)^2}.$
Substitution into equation \eqref{eq:block-scalar-equations2}, followed by elimination of denominators, gives
\begin{equation*}
 \eta\mu(1-\mu)(\widetilde{\mu}+1)^2
 +(\mu \eta^2-1)(\mu\widetilde{\mu}-1)(\widetilde{\mu}-\mu)=0.
\end{equation*}
Reorganizing by the powers of $\widetilde{\mu}$, we have
\begin{equation} \label{eq:determinant-polynomial}
 P_{\eta,\mu}(\widetilde{\mu})
 :=\big( \mu(\eta-1)(\eta\mu+1)\big)(\widetilde{\mu}^2+1)
 +(1+\mu^2+2\eta\mu-2\eta\mu^2-\eta^2\mu-\eta^2\mu^3)\widetilde{\mu}=0.
\end{equation}

\underline{Case 2, $\mu=-1$}. Equations~\eqref{eq:block-scalar-equations1} give $a=-d$ and $\widetilde{\mu}=-1$. Equation~\eqref{eq:block-scalar-equations2} then becomes
\begin{equation*}
 (\eta-1)^2d^2=\eta^2+1.
\end{equation*}
Since $\eta\neq1$, we may choose $d$ with $d^2=(\eta^2+1)/(\eta-1)^2$, and then set $a:=-d$. Meanwhile, $P_{\eta,-1}(\widetilde{\mu})=(\eta-1)^2(\widetilde{\mu}+1)^2$, so $\widetilde{\mu}=-1$ is also a solution of $P_{\eta,-1}(\widetilde{\mu})=0$.

Combining the two cases, the equation $X_2X_1^{-1}X_2^{-1}X_1=M_U$ reduces to $P_{\eta,\mu}(\widetilde{\mu})=0$. Under the hypotheses $\lambda^2\neq1, \lambda^2\mu\neq-1$, the quadratic and constant coefficients of $P_{\eta,\mu}(\widetilde{\mu})$ are nonzero. Hence, it is a genuine quadratic polynomial. 

Recall that the representation is diagonalizable $\murk$-one, which requires $\widetilde{\mu}\neq 1$.  We compute
\[
 P_{\eta,\mu}(1)
 =-(\mu-1)(\eta^2\mu^2-\eta^2\mu+4\eta\mu-\mu+1).
\]
Since $\mu\neq1$, the third hypothesis ensures that neither root is $1$. Finally, when $\mu\neq-1$, a root of $P_{\eta,\mu}(\widetilde{\mu})=0$ determines the desired values of $a,d$; when $\mu=-1$, we use the choice of $d$ made above and set $a=-d$.
\end{proof}

\begin{remark}
\label{rmk:doubling-hypotheses}
The three conditions in \eqref{eq:doubling-hypotheses} rule out the following cases:
\begin{enumerate}
	\item If $\lambda^2=1$, then $L_U =\pm I$, which gives $X_2=X_1$, and  $M_U = I$.
	\item If $\lambda^2 \mu= -1$ and $\lambda^2\neq1$, then the quadratic and constant terms of $P_{\eta,\mu}$ become zero, and the only solution is $\widetilde{\mu} =0$. In other words, $X_1$ is not invertible, and it cannot define a representation.
	\item If $ \lambda^4\mu^2-\lambda^4\mu+4\lambda^2\mu-\mu+1=0$, and $\lambda^2\neq1, \lambda^2\mu\neq-1$, then $\widetilde{\mu}=1$. The matrix $X_1$ is unipotent $\murk$-one but not diagonalizable $\murk$-one.
\end{enumerate}

\end{remark}

\begin{theorem}
\label{thm:one-step-doubling}
Let $\rho\colon\pi_K\to\GL(V)$ be an $n$-dimensional irreducible
diagonalizable $\murk$-one representation.  Suppose its peripheral data
$(\mu,\lambda)$ satisfy  \eqref{eq:doubling-hypotheses} in Lemma \ref{lem:whitehead-block-existence}.  Then there
exists a $2n$-dimensional irreducible diagonalizable
$\murk$-one representation
\begin{equation*}
 \widetilde\rho\colon\pi_{\Wh(K)}
 \longrightarrow\GL(V\oplus V^*).
\end{equation*}
\end{theorem}

\begin{proof}
Recall the decomposition $V\oplus V^*=U\oplus W$. The matrices $X_1, X_2$ from
Lemma~\ref{lem:whitehead-block-existence} act on $U$. Extend both by the
identity on $W$, and denote them by $\widetilde{X}_1, \widetilde{X}_2$. Define $\widetilde\rho$ on the Whitehead pattern by
\begin{equation*}
 \widetilde\rho(m_{\mathcal V}):=\rho^\#(m_K),\qquad
 \widetilde\rho(\ell_{\mathcal V}):=\rho^\#(\ell_K),\qquad
 \widetilde\rho(x_i):= \widetilde X_i\quad(i=1,2).
\end{equation*}
On the knot exterior, define $\widetilde\rho|_{\pi_K} = \rho^\#|_{\pi_K}$. By van Kampen, this is a well-defined  $2n$-dimensional representation of $\pi_{\Wh(K)}$.

To see the microlocal rank, we choose the meridian $x_1$ of $\Wh(K)$. By construction, $\widetilde{X}_1$ is the identity on $W$, and on $U$, it is straightforward to check that $e+e^*$ is an eigenvector of eigenvalue $\widetilde{\mu} \neq1$, and $(d-1)e+ (1-a) e^*$ is an eigenvector with the other eigenvalue, $1$. Therefore, $\widetilde\rho$ is a diagonalizable $\murk$-one representation.

Finally, it remains to prove irreducibility and identify the augmentation rank.
Since $\rho$ is irreducible, so is $\rho^*$.  Consequently, the two
nonzero invariant orbit spans satisfy
\begin{equation*}
 \Span\{\rho(g)e:g\in\pi_K\}=V,
 \qquad
 \Span\{\rho^*(g)e^*:g\in\pi_K\}=V^*.
 \label{eq:input-orbit-spanning}
\end{equation*}

Let $v_1: = e+e^*$ be the distinguished right $\widetilde{\mu}$-eigenvector of $\widetilde{X}_1$. Then the distinguished $\widetilde{\mu}$-eigenvector of $\widetilde{X}_2$ is $v_2:=\lambda^{-1}e+\lambda e^*$. In matrix form, they are column vectors $v_1 = (1,1)^t, v_2 = (\lambda^{-1}, \lambda)^t$. 
For any $g\in\pi_K$, the action of $\widetilde\rho(g)$ on the eigenvectors $v_1,v_2$ is 
\begin{equation*}
 \widetilde\rho(g)v_1=\rho(g)e\oplus\rho^*(g)e^*,
 \qquad
  \widetilde\rho(g)v_2=\lambda^{-1}\rho(g)e
          \oplus\lambda\rho^*(g)e^*.
\end{equation*}
Because $\lambda^2\neq1$, we have
\begin{equation*}
 \frac{\lambda \widetilde\rho(g)v_1-\widetilde\rho(g)v_2}{\lambda-\lambda^{-1}}
   =\rho(g)e\oplus0,
 \qquad
 \frac{\widetilde\rho(g)v_2-\lambda^{-1}\widetilde\rho(g)v_1}{\lambda-\lambda^{-1}}
   =0\oplus\rho^*(g)e^*.
\end{equation*}
Hence the right distinguished meridian orbit spans $V\oplus V^*$.

Similarly, regard the row vector $w_1:=(a-1,d-1)\in U^*$ as a
covector on $V\oplus V^*$ by extending it by zero on $W$.  Direct multiplication gives
$w_1X_1=((a-1)(a+d-1),(d-1)(a+d-1))=\widetilde{\mu} w_1$.  Since
$X_2=L_U^{-1}X_1L_U$, the row vector
$w_2:=w_1L_U=\bigl(\lambda(a-1),\lambda^{-1}(d-1)\bigr)$ satisfies
$w_2X_2=\widetilde{\mu} w_2$.

We claim that $a\neq1$ and $d\neq1$.  If $a=1$, then
\eqref{eq:block-scalar-equations1} gives $d=\mu^{-1}$, and
\eqref{eq:block-scalar-equations2} becomes
$\eta(1-\mu)\mu^{-1}=0$, a contradiction. Similarly, $d=1$ also leads to a contradiction. Thus both
coordinates of $w_1$ are nonzero.  Since $\lambda^2\neq1$, the same
separation calculation in
$(V\oplus V^*)^*=V^*\oplus V^{**}$, together with
\eqref{eq:input-orbit-spanning}, shows that the left distinguished
meridian orbit spans $(V\oplus V^*)^*$.

Since $w_1(v_1)=\widetilde{\mu}-1$, we normalize it to $v_1^* :=w_1/(\widetilde{\mu}-1)$, so as to pair
with $v_1$.  The preceding arguments imply
\begin{equation}
 \Span\{\widetilde\rho(h)v_1:h\in\pi_{\Wh(K)}\}=V\oplus V^*,
 \qquad
 \Span\!\left\{\widetilde\rho^*(h)v_1^*:
                  h\in\pi_{\Wh(K)}\right\}
   =(V\oplus V^*)^*.
 \label{eq:output-two-sided-spanning}
\end{equation}
Let $S\subset V\oplus V^*$ be $\widetilde\rho$-invariant.  If
$(\widetilde X_1-I)S\neq0$, then $S$ contains
$\operatorname{im}(\widetilde X_1-I)=\Bbbk v_1$, and the first equality in
\eqref{eq:output-two-sided-spanning} gives $S=V\oplus V^*$.  If
$(\widetilde X_1-I)S=0$, then $S$ lies in the $1$-eigenspace of
$\widetilde X_1$, which $v_1^*$ annihilates.  Invariance of $S$
implies that every
$\widetilde\rho^*(h)(v_1^*)$ also annihilates $S$, and the second equality in
 \eqref{eq:output-two-sided-spanning} gives $S=0$.  Thus
$\widetilde\rho$ is irreducible.
\end{proof}

\begin{corollary}
\label{cor:whitehead-rank-equality}
Let $K$ be a nontrivial knot admitting an $n$-dimensional irreducible
diagonalizable $\murk$-one representation whose peripheral data
$(\mu,\lambda)$ satisfy \eqref{eq:doubling-hypotheses} in Lemma \ref{lem:whitehead-block-existence}.
If $n=b(K)$, then
\[
 \operatorname{ar}(\Wh(K))
 =\operatorname{mr}(\Wh(K))
 =b(\Wh(K))=2n.
\]
\end{corollary}

\begin{proof}
By Theorem~\ref{thm:one-step-doubling}, there is a $2n$-dimensional
irreducible diagonalizable $\murk$-one representation $\widetilde\rho$
of $\pi_{\Wh(K)}$.  Use the eigenvectors $v_1$ and $w_1$ from the proof
of the theorem.  After enlarging the meridional generating set if
necessary, the matrix of the augmentation induced by $\widetilde\rho$
factors as
\begin{equation}
 R_{ij}=(1-\widetilde{\mu})
 \left(v_1^*\widetilde\rho(g_i)\right)
 \bigl(\widetilde\rho(g_j^{-1})v_1\bigr).
 \label{eq:doubled-augmentation-factorization}
\end{equation}
Equation~\eqref{eq:output-two-sided-spanning} says that the covectors and
vectors in \eqref{eq:doubled-augmentation-factorization} span
$(V\oplus V^*)^*$ and $V\oplus V^*$, respectively. Hence, the first
factor has column rank $2n$ and the second has row rank $2n$.  Since
$\widetilde{\mu}\neq1$, their product has rank $2n$.  Consequently,
\[
 \operatorname{ar}(\Wh(K))\geq2n.
\]

The Whitehead pattern has geometric index two.  A minimal bridge
presentation of $K$ gives a doubled bridge presentation of $\Wh(K)$,
while the satellite bridge inequality gives the reverse bound.  Since
$K$ is nontrivial, these two bounds apply and give
\[
 b(\Wh(K))=2b(K)=2n;
\]
see \cite[Section~4]{Dutra}.  Proposition~\ref{prop:rank-inequalities}
now yields
\[
 2n\leq\operatorname{ar}(\Wh(K))
 \leq\operatorname{mr}(\Wh(K))
 \leq b(\Wh(K))=2n,
\]
and all three ranks are equal to $2n$.
\end{proof}

\begin{proposition}\label{prop:induced-peripheral}
The peripheral data of $\widetilde\rho\colon\pi_{\Wh(K)}\to\GL(V\oplus V^*)$ are the following:
\begin{enumerate}
	\item If $\mu\neq-1$, then
	\begin{equation}
    \left(
    \widetilde \mu = a+d-1,\,
    \widetilde\lambda = \frac{\mu+\lambda^2-\widetilde{\mu}-\mu\lambda^2\widetilde{\mu}}
     {\mu\lambda^2-\mu\widetilde{\mu}-\lambda^2\widetilde{\mu}+1}
    \right).
    \label{eq:doubled-peripheral-generic}
    \end{equation}
    where $a, d$ are the choices made in the construction.
	
    \item If $\mu=-1$, then $a=-d$, and
	\begin{equation}
     \left(
     \widetilde\mu=-1,\,
     \widetilde\lambda
     =\eta+1+\eta^{-1}+(\eta-\eta^{-1})d
     \right),
     \label{eq:doubled-peripheral-exceptional}
    \end{equation}
    where $d$ is the choice made in the construction.

\end{enumerate}
\end{proposition}

\begin{proof}
The meridian of $\Wh(K)$ is $x_1$, and
$X_1v_1=\widetilde{\mu} v_1$.  By \eqref{eq:whitehead-longitude}, the action of the longitude on $U$ is
\begin{equation*}
 \widetilde\rho(\ell_{\Wh})|_U
 =X_1^{-1}L_UX_1L_U^{-1}X_1^{-1}L_U^{-1}X_1L_U.
\end{equation*}

Because $v_1 = e+e^*$ is the $\widetilde{\mu}$-eigenvector of $X_1$, it must be an eigenvector for $\widetilde\rho(\ell_{\Wh})$. A straightforward calculation yields
\begin{equation}
 \widetilde\rho(\ell_{\Wh})v_1
 =   \frac{\eta(a-1)+(d-1)}{(a-1)+\eta(d-1)} v_1.
  \label{eq:longitude-universal-formula}
\end{equation}
Note the denominator $(a-1)+\eta(d-1)$ is never zero. Otherwise, suppose 
$ (a-1)+\eta(d-1)=0$. 
Together with \(a=\mu d\), this gives
$
 (\mu+\eta)d=1+\eta.
$
If \(\mu+\eta=0\), then \(\eta=-1\) and hence \(\mu=1\), which is a
contradiction. Therefore,
$
 d=\frac{1+\eta}{\mu+\eta}.
$
Substituting this expression and $ a=\mu d $ into
\eqref{eq:block-scalar-equations2}, we obtain
\[
 0=\eta(1-\mu)ad+(\mu \eta^2-1)(a-1)(d-1)
  =\frac{\eta(1-\mu)(1+\mu \eta)^2}{(\mu+\eta)^2}.
\]
This yields either  \(\eta=0 \), \(\mu = 1\), or
\(\mu \eta =-1\), which is a contradiction. 

Next, we shall solve $a$ and $d$ in terms of $\mu$ and $\lambda$ (or $\eta$). There are two possibilities.

\underline{Case 1, $\mu\neq-1$}.  From $a=\mu d$ and
$\widetilde{\mu}=a+d-1$, we have
$a-1=(\mu\widetilde{\mu}-1)/(\mu+1)$ and
$d-1=(\widetilde{\mu}-\mu)/(\mu+1)$.  Substituting these identities into \eqref{eq:longitude-universal-formula} gives the desired
\eqref{eq:doubled-peripheral-generic}. 

\underline{Case 2, $\mu=-1$}.  Then $a=-d$, $\widetilde{\mu}=-1$, and
$(\eta-1)^2d^2=\eta^2+1$.  Substitution into
\eqref{eq:longitude-universal-formula} gives
\[
 \widetilde\lambda
 =\frac{\eta+1+(\eta-1)d}{\eta+1-(\eta-1)d}
 =\eta+1+\eta^{-1}+(\eta-\eta^{-1})d.
\]
This proves the desired 
\eqref{eq:doubled-peripheral-exceptional}.

\end{proof}

\subsection{Iteration and the algebraic peripheral correspondence}
\label{subsec:iteration}

The main conclusion of the preceding section is a rank-doubling
construction. Starting from an $n$-dimensional irreducible
diagonalizable $\murk$-one representation of $\pi_K$, under suitable conditions, we obtain a $2n$-dimensional
irreducible diagonalizable $\murk$-one representation of
$\pi_{\Wh(K)}$. If the construction can be repeated $N$ times, it
produces a $2^Nn$-dimensional representation of
$\pi_{\Wh^N(K)}$, and hence
$ \operatorname{ar}\bigl(\Wh^N(K)\bigr)\geq 2^Nn$.

% BEGIN INCLUDED ZL CONVENTIONS
Most of the constructions work for any field $\Bbbk$, but we eventually take $\Bbbk= \mathbb C$ for the geometric application. An affine algebraic set is a common zero locus of polynomials, with its Zariski topology; it need not be irreducible. We write $\Bbbk[A]$ for its coordinate ring. A regular function on an open subset is locally a quotient of polynomials with a nonvanishing denominator. An algebraic set is \emph{irreducible} if it is nonempty and not a union of two proper closed subsets. An irreducible component is a maximal irreducible closed subset. A nonempty open subset of an irreducible algebraic set is dense, and a morphism is \emph{dominant} if its image is dense. Dimension is the maximum length of a strict chain of nonempty irreducible closed subsets. A zero-dimensional irreducible algebraic set
is a point, and a curve is irreducible and one-dimensional. See
\cite[Chapter I, Sections 1--4]{hartshorne1977algebraic}.

A subset is \emph{constructible} if it is a finite union of locally closed
subsets. Chevalley's theorem says that the image of a constructible set
under a morphism is constructible
\cite[Chapter 1, Theorem 4.91 and Observation 4.92]{ferrersantos2005actions1}.
We shall also use the elementary consequence that a dense constructible
subset of an irreducible algebraic set contains a nonempty open subset.
These conventions apply to the peripheral tori and correspondences below.
% END INCLUDED ZL CONVENTIONS

Let $(\Bbbk^\times)^2$ be the algebraic torus with coordinates
$(\mu,\eta=\lambda^2)$.  Recall from Section~\ref{subsec:one-step-doubling} that
the Whitehead doubling construction requires $(\mu,\eta)\in\mathcal U$,
where
\begin{equation*}
 \mathcal U
 :=\left\{(\mu,\eta)\in(\Bbbk^\times)^2:
 (\mu-1)(\eta-1)(\mu \eta+1)
 (\mu^2\eta^2-\mu \eta^2+4\mu \eta-\mu+1)\neq0\right\}.
\end{equation*}
Here the first constraint is to ensure that $\rho(m_K)$ is diagonalizable, and the remaining three constraints are the hypotheses \eqref{eq:doubling-hypotheses} used in Lemma \ref{lem:whitehead-block-existence}, Theorem \ref{thm:one-step-doubling} and Corollary \ref{cor:whitehead-rank-equality}, with reasons explained in Remark \ref{rmk:doubling-hypotheses}. As we repeat the doubling construction,  the same conditions must hold for the output peripheral data at every stage. We begin with the definition for the loci where $N$ iterations are allowed. 

\begin{definition}
\label{def:admissible-doubling-chain}
Define a length-$N$ sequence consisting of
\[
 \big((\mu_j,\eta_j)_{0\leq j\leq N}, (a_j,d_j)_{0\leq j<N}\big)
 \in(\Bbbk^\times)^{2N+2}\times \Bbbk^{2N},
\]
satisfying the following conditions:

\begin{enumerate}
	\item doubling constraints: for $0\leq j<N$,
	\[
	 \mu_j\neq1,\qquad \eta_j\neq1,
	 \qquad \eta_j\mu_j\neq-1, \qquad \eta_j^2\mu_j^2-\eta_j^2\mu_j+4\eta_j\mu_j-\mu_j+1\neq 0;
	\]
	\item internal relations: for $0\leq j<N$,
	\[
	 a_j =\mu_jd_j, \qquad
	 \eta_j(1-\mu_j)a_jd_j
	   +(\mu_j \eta_j^2-1)(a_j-1)(d_j-1) =0.
	\]
	
	\item inductive relations: for $0\leq j<N$,
	\[
 \mu_{j+1}:=a_j+d_j-1,
 \qquad
 \eta_{j+1}:=\left(
   \frac{\eta_j(a_j-1)+(d_j-1)}{(a_j-1)+\eta_j(d_j-1)}
 \right)^2.
	 \]
\end{enumerate}
\end{definition}

Let $\mathcal{C}_N$ be the set of length-$N$ sequences. Define the projection
\begin{align*}
 \operatorname{pr}:\quad
(\Bbbk^\times)^{2N+2}\times \Bbbk^{2N} &\longrightarrow(\Bbbk^\times)^2,\\
\left((\mu_j,\eta_j), (a_j,d_j)\right) &\longmapsto(\mu_0,\eta_0).
\end{align*}
Define $\mathcal S_N:=\operatorname{pr}(\mathcal C_N)$ to be the image. Note that
$\mathcal S_0=(\Bbbk^\times)^2$ and $\mathcal S_1=\mathcal U$.
Moreover, by truncating a length-$(N+1)$ sequence, we have
\begin{equation*}
 \mathcal S_{N+1}\subseteq\mathcal S_N.
\end{equation*}

A length-$N$ sequence records the compatible data required to carry out
$N$ successive doubling steps. Suitable initial peripheral data would start the iteration.

\begin{proposition}
\label{thm:finite-depth-iteration}
Let $\rho_0\colon\pi_K\to\GL(V_0)$ be an $n$-dimensional irreducible
diagonalizable $\murk$-one representation with peripheral data
$(\mu_0,\lambda_0)$. Suppose $(\mu_0,\lambda_0^2)\in\mathcal S_N$. Then for each $0\leq j\leq N$, there exists a $2^jn$-dimensional irreducible diagonalizable $\murk$-one representation $\rho_j\colon\pi_{\Wh^j(K)}\longrightarrow\GL(V_j)$.
\end{proposition}

\begin{proof}
Choose a length-$N$ sequence $\big((\mu_j,\eta_j)_{0\leq j\leq N}, (a_j,d_j)_{0\leq j<N}\big)$ representing
$(\mu_0,\lambda_0^2)\in\mathcal S_N$. For $0\leq j<N$, use the block determined by $a_j,d_j$ in the
construction of Theorem~\ref{thm:one-step-doubling}, and define
$V_{j+1}:=V_j\oplus V_j^*$.  The doubling constraints allow the
construction at every stage, while the inductive relations and the
peripheral formula give
$\mu_{j+1}=a_j+d_j-1$ and
$\lambda_{j+1}^2=\eta_{j+1}$.  Induction therefore gives irreducible
diagonalizable $\murk$-one representations
\begin{equation*}
 \rho_j\colon\pi_{\Wh^j(K)}\longrightarrow\GL(V_j),
 \qquad 0\leq j\leq N.
\end{equation*}
Since $\dim V_{j+1}=2\dim V_j$, we have $\dim V_j=2^jn$.
\end{proof}

\begin{corollary}
\label{cor:finite-depth-rank-equality}
Let $K$ be nontrivial, and let
$\rho_0\colon\pi_K\to\GL(V_0)$ be an $n$-dimensional irreducible diagonalizable $\murk$-one representation with peripheral data $(\mu_0,\lambda_0)$. If $b(K)=n$ and $(\mu_0,\lambda_0^2)\in\mathcal S_N$, then
\begin{equation}
 \operatorname{ar}\bigl(\Wh^j(K)\bigr)
 =\operatorname{mr}\bigl(\Wh^j(K)\bigr)
 =b\bigl(\Wh^j(K)\bigr)
 =2^j n,
 \qquad 0\leq j\leq N.
\end{equation}
\end{corollary}

\begin{proof}
This proof iterates the argument in Corollary \ref{cor:whitehead-rank-equality} $N$ times. Namely, Proposition~\ref{thm:finite-depth-iteration} yields
$
 \operatorname{ar}\bigl(\Wh^j(K)\bigr)\geq2^j n.
$
Iterating the doubled bridge presentation gives
$
 b\bigl(\Wh^j(K)\bigr)\leq2^j b(K)=2^j n.
$
Together with Proposition~\ref{prop:rank-inequalities}, these bounds
give
\[
 2^j n
 \leq\operatorname{ar}\bigl(\Wh^j(K)\bigr)
 \leq\operatorname{mr}\bigl(\Wh^j(K)\bigr)
 \leq b\bigl(\Wh^j(K)\bigr)
 \leq2^j n,
\]
which proves the claimed equalities.
\end{proof}

\begin{example}
We now consider the case when $K = 3_1$ is the right-handed trefoil, which has bridge index $2$. A routine computation shows that we can express its knot group as
\[
    \pi_K = \langle a,b~|~aba = bab \rangle.
\]
Furthermore, we can identify $m_K = a$ and $\ell_K = (ab)^3a^{-6}$. We can write down a concrete representation
\[
    \rho: \pi_K \to {\rm GL}(\mathbb{C}^2)
\]
as
\[
    \rho(a) = \left(
		\begin{matrix}
			2 & 1 \\
			0 & 1
		\end{matrix}
	\right) \aand \rho(b) = \left(
		\begin{matrix}
			1 & 0 \\
			-2 & 2
		\end{matrix}
	\right)
\]
One can verify it is indeed a diagonalizable $\mu{\rm rk}$-one representation with peripheral data $(\mu_0,\lambda_0) = (2,-\frac{1}{8})$. At each iteration, there are two choices of $\mu$, which are roots of the quadratic polynomial in \eqref{eq:determinant-polynomial}. Since the two distinct candidates are inverse to each other, we always take the one whose norm is greater than $1$. This provides an iteration $(\mu_n,\lambda_n^2)\to (\mu_{n+1},\lambda_{n+1}^2)$ provided $(\mu_n,\lambda_n^2)\in\mathcal{S}_1$. A short Python script prints the peripheral data of the first five iterations as follows.
\begin{table}[ht]
%\color{orange}
\centering
\begin{tabular}{rcc}
\hline
$n$ & $\mu_n$ & $\lambda_n$ \\
\hline
1 & $1.906106187303503$ & $-1.777215268543279\times 10^{-2}$ \\
2 & $1.904356477557643$ & $-3.494634031728678\times 10^{-4}$ \\
3 & $1.904355802094616$ & $-1.350405429773338\times 10^{-7}$ \\
4 & $1.904355802094515$ & $-2.016457262216791\times 10^{-14}$ \\
5 & $1.904355802094515$ & $-4.496128493818061\times 10^{-28}$ \\
\hline
\end{tabular}
\par\vspace{6pt}
\caption{Peripheral data for first 5 iterations, choosing $\mu_n>1$ at every step. Values are rounded to 16 significant digits.}
\label{tab:trefoil-whitehead-iterations}
\end{table}
With some effort, one could show that the initial data $(2,\frac{1}{64})$ lies in $\mathcal{S}_N$ for all $N$. The final remark is to note that this infinite iteration stays ``bounded'' in contrast to the iterable branches discussed in Section \ref{sec:escaping-curves}.
\end{example}

Note that the above discussion does not involve the consideration of a specific knot $K$. The sets $\mathcal S_N$ depend only on the Whitehead pattern, not on the companion knot. To apply the construction to a given companion $K$ for a prescribed number $N$ of steps, one must find a representation of the required dimension whose peripheral data satisfy $(\mu,\lambda^2)\in \mathcal S_N$. We next show that each $\mathcal S_N$ contains a Zariski-open dense
subset of the ambient peripheral torus.  This statement does
not by itself imply that $\mathcal S_N$ meets a prescribed peripheral
curve of $K$. The next proposition shows an essential property of $\mathcal S_N$, which asserts enough general points. 

\begin{proposition}
\label{prop:finite-depth-chain-loci}
For every $N\geq0$, the set $\mathcal S_N$ is constructible and contains a Zariski-open dense subset of $(\Bbbk^\times)^2$. Moreover, since $\Bbbk =\mathbb{C}$ is uncountable, $\cap_{N\in \mathbb{N}}\mathcal S_N$ is Zariski dense.
\end{proposition}

To prove this proposition, we first consider the correspondence defined by the doubling construction. Consider two tori $(\Bbbk^\times)^2_{({\mu}, {\eta})}$ and $(\Bbbk^\times)^2_{(\widetilde{\mu}, \widetilde {\lambda})}$. Let $\pi_1, \pi_2$ be the projections from the product to the first and second components.  We also define ${\widetilde{\pi}_2}$ by
\[
{\widetilde{\pi}_2}(\mu,\eta,\widetilde{\mu}, \widetilde {\lambda}):=(\widetilde{\mu}, \widetilde {\lambda}^2)
\]
These maps fit into the commutative diagram

\begin{center}
\begin{tikzpicture}[
  obj/.style={align=center,inner sep=2pt,font=\small},
  map/.style={-{Latex},semithick}
]
 \node[obj] (input) at (-2,-0.3)
   {$(\Bbbk^\times)^2_{(\mu,\eta)}$};
 \node[obj] (W) at (0,1.2) {$(\Bbbk^\times)^2_{({\mu}, {\eta})}\times (\Bbbk^\times)^2_{(\widetilde{\mu}, \widetilde {\lambda})}$};
 \node[obj] (output) at (1.5,-0.3)
    {$(\Bbbk^\times)^2_{(\widetilde{\mu}, \widetilde {\lambda})}$};
 \node[obj] (next) at (4.5,-0.3)
   {{$(\Bbbk^\times)^2_{(\widetilde{\mu}, \widetilde {\eta})}$}};
 \draw[map] (W) -- node[above] {$\scriptstyle{\pi_1}$} (input);
 \draw[map] (W) -- node[above] {$\scriptstyle{\pi_2}$} (output);
 \draw[map] (W) -- node[above] {$\scriptstyle{\widetilde{\pi}_2}$} (next);
 \draw[map] (output) --
   node[below] {$\scriptstyle{\widetilde {\eta} = \widetilde {\lambda}^2}$}
   (next);
\end{tikzpicture}
\end{center}

We formulate a correspondence between $(\mu, \eta)$ and $(\widetilde\mu,\widetilde\lambda)$ from the iteration procedure. Recall from \eqref{eq:determinant-polynomial} that $\widetilde\mu$ is the solution of $P_{\eta,\mu}(\widetilde\mu) =0$, and from \eqref{eq:doubled-peripheral-generic} that $\lambda$ is a rational polynomial. Let $\delta(\mu,\eta,\widetilde{\mu}):=\mu \eta-\mu\widetilde{\mu}-\eta\widetilde{\mu}+1$ be the denominator, then we can rewrite the expression as a polynomial $\delta\widetilde\lambda +\mu \eta\widetilde\mu-\mu-\eta+\widetilde\mu =0$. Then we define:
\begin{equation}
 \widehat{\mathcal W}
 :=\left\{(\mu,\eta,\widetilde\mu,\widetilde\lambda)
 :P_{\eta,\mu}(\widetilde\mu) =0,\, 
 \delta\widetilde\lambda +\mu \eta\widetilde\mu-\mu-\eta+\widetilde\mu =0
 \right\}
 \subset (\Bbbk^\times)^2_{(\mu,\eta)}
 \times(\Bbbk^\times)^2_{(\widetilde\mu,\widetilde\lambda)}.
 \label{eq:raw-whitehead-correspondence}
\end{equation}
Note that $P_{\eta,\mu}(\widetilde\mu) =0$ is generically quadratic, namely
\[
P_{\eta,\mu}(\widetilde\mu) = A\widetilde{\mu}^2 + B\widetilde{\mu} + A,
\]
where
\begin{align*}
    A &= \mu(\eta-1)(\eta\mu+1), \\
    B &= 1+\mu^2+2\eta\mu-2\eta\mu^2-\eta^2\mu-\eta^2\mu^3.
\end{align*}
The discriminant is 
\begin{equation}
    \label{eq:peripheral-quadratic-discriminant}
B^2 - 4A^2 = (\mu-1)(\mu+1)^2(\mu \eta^2-1)
  (\mu^2\eta^2-\mu \eta^2+4\mu \eta-\mu+1).
\end{equation}
Under the conditions $\mu \neq -1$ and $\mu \eta^2 \neq 1$, the discriminant is non-zero. Hence, we introduce the Zariski-open subset
\begin{equation*}
 \mathcal U':=\bigl\{(\mu,\eta)\in \mathcal U: \mu \neq -1,\,  \mu \eta^2 \neq 1\bigr\}\subset \mathcal U \subset(\Bbbk^\times)^2_{(\mu,\eta)}.
\end{equation*}
Finally, we introduce
\[ 
\mathcal W := \text{Closure}(\widehat{\mathcal W}\cap \pi_1^{-1}(\mathcal U'))\cap \widehat{\mathcal W}.
\]
This is the correspondence that we will work with.

\begin{proposition}
\label{prop:peripheral-correspondence-dominance}
We have the following properties:
\begin{enumerate}
    \item $\mathcal W$ is irreducible.
	\item $\pi_1|_{\mathcal W}: \mathcal W \rightarrow (\Bbbk^\times)^2_{(\mu,\eta)}$ is dominant and generically two-to-one.
	\item $\pi_2|_{\mathcal W}:\mathcal W\to (\Bbbk^\times)^2_{(\widetilde\mu,\widetilde\lambda)}$ is dominant and generically at most two-to-one.
	\item $\widetilde{\pi}_2|_{\mathcal W}:\mathcal W\to (\Bbbk^\times)^2_{(\widetilde\mu,\widetilde \eta)}$ is also dominant.
\end{enumerate}

\end{proposition}

\begin{proof}
(1) Over the Zariski-open dense subset $\mathcal U'$, $P_{\eta,\mu}(\widetilde\mu) =0$ has two distinct solutions. Moreover, the discriminant $B^2-4A^2$ is not a square in $\Bbbk(\mu,\eta)$, and $P_{\eta,\mu}$ is irreducible over $\Bbbk(\mu,\eta)$. Indeed, to see the non-squareness, the factor $(\mu-1)$ has order 1. The other factors $(\mu+1)$, $(\mu \eta^2-1)$, $(\mu^2\eta^2-\mu \eta^2+4\mu \eta-\mu+1)$ do not contain any power of $(\mu-1)$, because when evaluated at $\mu =1$, they equal $2, (\eta^2-1), 4\eta$ respectively, none of which is zero.

The coefficients of $P_{\eta,\mu}(\widetilde\mu)$ have no common factors, making it primitive. By Gauss's lemma, $P_{\eta,\mu}(\widetilde\mu)$ is also irreducible over $\Bbbk[\mu,\eta]$. Once $\eta, \mu, \widetilde\mu$ are fixed, $\widetilde\lambda$ is uniquely determined as a rational polynomial. In particular, over the Zariski-open dense locus ${\delta\neq 0}$, the variety is irreducible. The exceptional case, since we transformed the equation of $\widetilde\lambda$ into a polynomial, is that both the numerator and denominator equal zero, i.e. $\mu \eta\widetilde\mu-\mu-\eta+\widetilde\mu =0, \mu \eta-\mu\widetilde{\mu}-\eta\widetilde{\mu}+1=0.$ The solutions are $\mu = \widetilde\mu =1$ and $\eta$ can be anything. However, this is excluded from $\mathcal U$, and certainly excluded from $\mathcal U'$. Altogether, we proved that $\mathcal W$ is irreducible.

(2) Since $\mathcal U'$ is open dense in $(\Bbbk^\times)^2_{(\mu,\eta)}$, $\pi_1|_{\mathcal W}$ is dominant. A generic fiber has two choices of $\widetilde\mu$ and a corresponding $\widetilde\lambda$.

(3) To study the property of $\pi_2$, we need to express $\eta,\mu$ in terms of $\widetilde\mu, \widetilde\lambda$. First, with $\mu, \widetilde\mu, \widetilde\lambda$ fixed, $\eta$ is uniquely determined by
\begin{equation*}
 \eta=\frac{\mu\widetilde\mu\widetilde\lambda
          +\mu-\widetilde\mu-\widetilde\lambda}
         {\mu\widetilde\mu+\mu\widetilde\lambda
          -\widetilde\mu\widetilde\lambda-1}.
\end{equation*}
Substituting in $P_{\eta,\mu}(\widetilde\mu)=0$ and clearing denominators, we get $(\mu-1)(\mu+1)^2 Q_{\widetilde\mu,\widetilde\lambda}(\mu)=0$, where
\begin{equation*}
Q_{\widetilde\mu,\widetilde\lambda}(\mu)
=\widetilde\mu(\widetilde\mu\widetilde\lambda+1)^2(\mu^2+1)
+(-\widetilde\mu^4\widetilde\lambda^2
+\widetilde\mu^4\widetilde\lambda
-2\widetilde\mu^3\widetilde\lambda
-\widetilde\mu^2\widetilde\lambda^2
-2\widetilde\mu^2\widetilde\lambda
-\widetilde\mu^2
-2\widetilde\mu\widetilde\lambda
+\widetilde\lambda
-1)\mu 
\end{equation*}
For generic values of $(\widetilde\mu,\widetilde\lambda)$, this is a quadratic polynomial. Moreover, since the quadratic and constant coefficients are equal, the two solutions multiply to $1$. In particular, they are inverse to each other and non-zero. Therefore, $\pi_2|_{\mathcal W}:\mathcal W\to (\Bbbk^\times)^2_{(\widetilde\mu,\widetilde\lambda)}$ is dominant
 and generically at most two-to-one.

(4) Since $\widetilde{\pi}_2$ is the composite of $\pi_2$ with the finite dominant map $(\widetilde\mu,\widetilde\lambda)\mapsto
(\widetilde\mu,\widetilde\lambda^2)$, $\widetilde{\pi}_2$ is also dominant.
\end{proof}

We next identify the part of the correspondence on which the one-step
construction is nondegenerate. Define 
\begin{equation}
 \mathcal W_{\mathrm{step}}
 :=\left\{w=(\mu,\eta, \widetilde\mu, \widetilde\lambda)\in\mathcal W:
 \pi_1(w)\in \mathcal U,\
 \delta(\mu,\eta, \widetilde\mu)\neq0\right\}.
 \label{eq:good-whitehead-correspondence}
\end{equation}
Observe that over $\mathcal{W}_{\mathrm{step}} $, we have $\mu\neq -1$. 
Indeed, if $\mu=-1$, then $\eta\neq1$ and the equation
\[
 P_{\eta,-1}(\widetilde\mu)
 =(\eta-1)^2(\widetilde\mu+1)^2
\]
forces $\widetilde\mu=-1$.  But then
$\delta(-1,\eta,-1)=0$, contrary to the definition of
$\mathcal W_{\mathrm{step}}$.

We also note that over $\mathcal{W}_{\mathrm{step}} $, given the input peripheral data $(\mu,\eta)$, the pairs $(a,d)$ and $(\widetilde\mu,\widetilde\lambda)$ determine each other in a unique way. Namely, by Proposition \ref{prop:induced-peripheral}, we have
\begin{equation*}
\widetilde\mu = a+d-1, \qquad
\widetilde\lambda = \frac{\eta(a-1)+(d-1)}{(a-1)+\eta(d-1)},
 \end{equation*}
 And conversely, we can solve $(a,d)$ as
\begin{equation*}
 a=\frac{\mu(\widetilde\mu+1)}{\mu+1}   \qquad
 d=\frac{\widetilde\mu+1}{\mu+1}.
\end{equation*}

We now iterate the correspondence with the above mild restrictions. Define
\begin{equation*}
 \mathcal G_0:=(\Bbbk^\times)^2,    \qquad
 \mathcal G_{N+1}
 :=\pi_1\!\left(
     \mathcal W_{\mathrm{step}}\cap \widetilde\pi_2^{-1}(\mathcal G_N)
   \right).
\end{equation*}

It follows from the inductive construction that a point $(\mu_0,\eta_0)$ lies in $\mathcal G_N$ if and only if there exist $w_j\in\mathcal W_{\mathrm{step}}$ for $0\leq j<N$ and $\{(\mu_j,\eta_j)\}_{1\leq j\leq N}$ such that
\begin{equation*}
 \pi_1(w_j)=(\mu_j,\eta_j),
 \qquad
 \widetilde\pi_2(w_j)=(\mu_{j+1},\eta_{j+1}),
 \qquad 0\leq j<N.
\end{equation*}
Therefore, by truncation, we find that 
\[
\mathcal G_{N+1}\subseteq\mathcal G_N.
\]
Moreover, given the data of $(\mu_j, \eta_j)$, we can solve $(a_j,d_j)$ from  $(\mu_{j+1}, \eta_{j+1})$. Hence, we can inductively define the length-$N$-sequence. This yields 
\begin{equation}\label{eq: G inside S}
	\mathcal G_N\subseteq\mathcal S_N.
\end{equation}

\begin{proposition}
\label{prop:generic-finite-depth-chains}
For every $N\geq0$, the set $\mathcal G_N$ is constructible and
contains a Zariski-open dense subset of $(\Bbbk^\times)^2$. 
\end{proposition}

\begin{proof}
Since $\widehat{\mathcal W}\cap \pi_1^{-1}(\mathcal U')\subset\mathcal W_{\mathrm{step}}$, we have that $\mathcal W_{\mathrm{step}}$ is a Zariski-open dense subset of $\mathcal W$. We prove by induction. The base case $N=0$ is clear. Suppose $\mathcal G_N$ is
constructible and contains a dense open subset
$U_N\subset(\Bbbk^\times)^2$.  The set
$\mathcal W_{\mathrm{step}}\cap\widetilde\pi_2^{-1}(\mathcal G_N)$ is
constructible, and Chevalley's theorem further shows that its image under $\pi_1$, which is $\mathcal G_{N+1}$, is also constructible. By Proposition~\ref{prop:peripheral-correspondence-dominance}, $\tilde{\pi}_2$ is dominant and $\mathcal W$ is irreducible. Therefore, $ \mathcal W_{\mathrm{step}}\cap\widetilde\pi_2^{-1}(U_n)\subset \mathcal W$ is Zariski-open dense. Since $\pi_1$ is dominant, its image under $\pi_1$ is also Zariski-open dense. Therefore, $\mathcal G_{N+1}$ is dense and constructible, which completes the induction.
\end{proof}

\begin{proof}[Proof of Proposition \ref{prop:finite-depth-chain-loci}]
Following Definition~\ref{def:admissible-doubling-chain}, $\mathcal C_N$ is defined by a finite collection of polynomials and non-vanishing conditions, hence it is locally closed and further constructible. By the Chevalley's theorem, its projection $\mathcal S_N = \operatorname{pr}(\mathcal C_N)$ is also constructible. By Proposition~\ref{prop:generic-finite-depth-chains}, the set
$\mathcal G_N$ contains a Zariski-open dense subset. Since $\mathcal G_N\subset \mathcal S_N$, $\mathcal S_N$ also contains a Zariski-open dense subset.

Finally, we show that $\cap_{N\in \mathbb{N}}\mathcal S_N$ is Zariski dense. For each $N$, we choose a non-zero polynomial $f_N$ such that its principal open subset $U_N :=\{f_N\neq 0\}\cap (\Bbbk^\times)^2\subset S_N$. We claim that $\cap_N U_N$ is Zariski dense. It suffices to show that for any principal open subset $\{g\neq 0\}$, the intersection $\{f_N\neq 0\}\cap (\cap U_N)\neq \emptyset$. Consider any element $h\in \{g, f_0,f_1,\dotsb\}$. Let $Z(h)\subset \Bbbk^\times_\mu$ be the subset of $\mu_0$ such that $h(\mu_0, \eta)=0$ with respect to $\eta$. Then $Z(h)$ is finite, since for $h = \sum_{i\geq 0} h_i(\mu)\eta^i$, each $h_i$ has finitely many roots. Note that $Z(g)\cup (\cup_N Z(f_N))$ is countable. Since $\Bbbk = \mathbb{C}$ is uncountable, we can choose $\mu_0\in \Bbbk^\times \setminus (Z(g)\cup (\cup_N Z(f_N)))$, such that $h(\mu_0, \eta)$ is non-zero. Iterate the argument for the second variable, we find $\eta_0$ such that neither $g(\mu_0,\eta_0)$ nor any $f_N(\mu_0,\eta_0)$ vanished. This completes the proof.
\end{proof}

% END INCLUDED SOURCE: sections/02_doubling.tex
% BEGIN INCLUDED SOURCE: sections/03_meridional_bounds.tex
\subsection{Meridional bounds}\label{sec:meridional-bounds}
We next compare the representations constructed in Section~\ref{subsec:iteration} with meridional generating sets. The upper bound below does not require an augmentation or a condition on its peripheral data.

\begin{proposition}\label{prop:mr-doubles-upper}
For every knot $K$ and every integer $N\geq0$,
\[
 \mr(\Wh(K))\leq2\mr(K),\qquad
 \mr(\Wh^N(K))\leq2^N\mr(K).
\]
\end{proposition}

\begin{proof}
Let $r=\mr(K)$. Choose meridional generators
$m_i=h_im_Kh_i^{-1}$, $1\leq i\leq r$, with $h_1=1$.
In the satellite presentation of Proposition~\ref{prop:van-kampen}, we have
\[
 m_{\mathcal V}=x_2x_1^{-1}x_2^{-1}x_1=m_K,
 \qquad \ell_{\mathcal V}=\ell_K,
 \qquad \ell_{\mathcal V}x_2\ell_{\mathcal V}^{-1}=x_1.
\]
In particular, both $x_1$ and $x_2$ are meridians of $\Wh(K)$.
Let $H$ be the subgroup generated by the $2r$ meridians
\[
 u_i=h_ix_1h_i^{-1},\qquad v_i=h_ix_2h_i^{-1},
 \qquad 1\leq i\leq r.
\]
Here elements of the companion group are understood through their images in
the satellite group. For each $i$,
\[
 v_iu_i^{-1}v_i^{-1}u_i
 =h_i(x_2x_1^{-1}x_2^{-1}x_1)h_i^{-1}=m_i.
\]
Thus $H$ contains the image of the companion group, and hence $\ell_{\mathcal V}$.
It also contains $x_1=u_1$ and $x_2=v_1$. These elements generate the
pattern group. The satellite presentation therefore gives
$H=\pi_{\Wh(K)}$, proving the first inequality. Induction gives the second.
The same presentation argument applies when $K$ is the unknot; it does not
require the companion homomorphism to be injective.
\end{proof}

Together with Proposition~\ref{prop:rank-inequalities}, this gives the unconditional bounds
\begin{equation}\label{eq:unconditional-meridional-bounds}
 \ar(\Wh^N(K))\leq\mr(\Wh^N(K))\leq2^N\mr(K).
\end{equation}
The next result gives a sufficient condition for equality. Recall the chain-admitting set $\mathcal S_N$ from
Definition~\ref{def:admissible-doubling-chain} and the paragraph following it.

\begin{theorem}\label{thm:mr-version}
Let $K$ be a knot and set $r=\mr(K)$. Suppose that $K$ admits a rank-$r$
augmentation $\epsilon$ with $\epsilon(\mu)\neq1$ and
\[
 \bigl(\epsilon(\mu),\epsilon(\lambda)^2\bigr)\in\mathcal S_N.
\]
Then
\[
 \ar(\Wh^N(K))=\mr(\Wh^N(K))=2^N\mr(K).
\]
\end{theorem}

\begin{proof}
The augmentation--representation correspondence recalled in
Section~\ref{subsec:augmentation-representations} gives an irreducible
$r$-dimensional diagonalizable $\murk$-one representation with the prescribed peripheral data.
Proposition~\ref{thm:finite-depth-iteration} then gives
$\ar(\Wh^N(K))\geq2^Nr$. Combining this with
\eqref{eq:unconditional-meridional-bounds} proves the result.
\end{proof}

\begin{remark}
The condition on the peripheral point is essential. In
Section~\ref{sec:first-doubles}, we show that the square knot $S$ satisfies
$\ar(S)=\mr(S)=b(S)=3$, whereas $\ar(\Wh(S))\leq5$.
Thus maximal augmentation rank for the companion alone does not suffice.
\end{remark}
% END INCLUDED SOURCE: sections/03_meridional_bounds.tex
% BEGIN INCLUDED SOURCE: sections/04_ZL_preliminaries.tex
\section{Algebraic curves and character varieties}\label{sec:algebraic-preliminaries}

Theorem~\ref{thm:mr-version} reduces the meridional-rank problem to
finding an augmentation whose peripheral point lies in $\mathcal S_N$.
To make this choice for every $N$, we use a curve of peripheral
values. We first recall how a point at infinity
of such a curve determines the orders of poles and zeros. We then record
the character-variety constructions used to obtain the curve from a
boundary slope. Throughout this section, the ground field is $\mathbb C$,
and we use the algebraic conventions of Section~\ref{subsec:iteration}.

\subsection{Curves, ideal points, and valuations}\label{subsec:curve-preliminaries}

For an irreducible algebraic set $A$, we denote its field of rational
functions by $\mathbb C(A)$. A rational map is defined on a nonempty open
subset, not necessarily on its whole source; two representatives define
the same map if they agree on their
common domain. A birational map restricts to an isomorphism between
nonempty open subsets. We write $\operatorname{Dom}(\varphi)$ for the
largest open subset on which a rational map $\varphi$ is defined.

\begin{definition}\label{defn: smooth projective model}
	Suppose $C$ is an affine algebraic curve. A \emph{smooth projective model} $(\widetilde{C},\tilde{\phi})$ of $C$ is a smooth projective curve $\widetilde{C}$ together with a birational map $\tilde{\phi}:\widetilde{C}\dashrightarrow C$. A point $P\in \widetilde{C}$ not in the domain of $\tilde{\phi}$ is called an ideal point of $C$.
\end{definition}

Every affine algebraic curve has a smooth projective model
\cite[Chapter I, Corollary 6.11]{hartshorne1977algebraic}.
The birational map identifies $\mathbb C(C)$ with
$\mathbb C(\widetilde C)$; we use this identification for rational
functions throughout. The following lemma is important to lift valuation-related results from the $\PSL_2(\mathbb{C})$ setup to the $\SL_2(\mathbb{C})$ setup.
\begin{lemma}\label{lem: surjective morphism on smooth projective model}
	Suppose $f: C\to D$ is a morphism between two affine algebraic curves so that $D = \overline{\im(f)}$. Suppose $(\widetilde{C},\tilde{\phi}_C)$ and $(\widetilde{D},\tilde{\phi}_D)$ are the respective smooth projective models. Then there exists a unique surjective morphism $\tilde{f}:\widetilde{C}\to \widetilde{D}$ such that as rational maps we have $\tilde{\phi}_D\circ \tilde{f} = f\circ \tilde{\phi}_C.$
\end{lemma}
\begin{proof}
	The existence of a dominant morphism $\tilde{f}: \widetilde{C}\to\widetilde{D}$ with the identity $\tilde{\phi}_D\circ \tilde{f} = f\circ \tilde{\phi}_C$ is established in \cite[Chapter I, Corollary 6.12]{hartshorne1977algebraic}. The meaning of ``dominant'' is precisely what we assumed ($D = \overline{\im(f)}$). \cite[Chapter 1, Theorem 4.100]{ferrersantos2005actions1} says that $\widetilde{C}$ is complete, and hence by \cite[Chapter 1, Theorem 4.98]{ferrersantos2005actions1} the image of $\tilde{f}$ is closed, and must be all of $\widetilde{D}$.
\end{proof}

Let $P\in\widetilde C$. Denote by $\mathcal O_{\widetilde C,P}$ the
local ring of rational functions regular near $P$, and by
$\mathfrak m_{\widetilde C,P}$ its maximal ideal of functions vanishing
at $P$. Since $\widetilde C$ is a smooth curve, this local ring is a
discrete valuation ring
\cite[Chapter I, Theorems 5.1 and 6.2A]{hartshorne1977algebraic}.
In particular, $\mathfrak m_{\widetilde C,P}$ has a generator $s$,
and every nonzero $f\in\mathbb C(\widetilde C)$ has a unique expression
$f=s^k h$ with $k\in\mathbb Z$ and $h$ a unit in
$\mathcal O_{\widetilde C,P}$. The integer $k$ is independent of the
choice of $s$.

\begin{definition}\label{defn: local parameter}\label{defn: valuation}
A generator $s$ of $\mathfrak m_{\widetilde C,P}$ is a \emph{local
parameter} at $P$. For $f=s^k h$ as above, its \emph{valuation} at $P$
is $\nu_P(f)=k$, and its \emph{pole order} is
\[
 \operatorname{pole}_P(f)=\max\{0,-\nu_P(f)\}.
\]
We set $\nu_P(0)=+\infty$. Thus zeros have positive valuation and poles
have negative valuation.
\end{definition}

The next lemma explains why the inequalities between valuations used
later are preserved when we pass to a covering curve.

\begin{lemma}\label{lem: valuation and surjective morphism}
	Suppose $\psi: \widetilde{C}\to \widetilde{C}'$ is a surjective morphism between two smooth projective curves. Let $P\in \widetilde{C}$ and $Q = \psi(P)$. Let $s_Q$ be the local parameter at $Q$ as in Definition~\ref{defn: local parameter}. Write
	\[
		e(P/Q) = \nu_P(\psi^*(s_Q)).
	\]
	Then the following statements are true:
	\begin{itemize}
		\item The value $e(P/Q)$ is a positive integer and independent of the choice of the local parameter $s_Q$.
		\item For any nonzero rational function $f$ on $\widetilde{C}'$, we have
		\[
			\nu_P(\psi^*f)=e(P/Q)\nu_Q(f).
		\]
		Consequently, we have
		\[
			\operatorname{pole}_P(\psi^*f)=e(P/Q)\operatorname{pole}_Q(f).
		\]
	\end{itemize}
\end{lemma}
\bpf
Since by construction $s_Q(Q) = 0$, we have $(\psi^*s_Q)(P) = 0$. \cite[Chapter I, Theorem 4.4]{hartshorne1977algebraic} implies that $\psi^*s_Q$ is not identically $0$, and hence $e(P/Q) = \nu_P(\psi^*(s_Q))$ is a positive integer. It is clearly independent of the choice of $s_Q$: two different choices differ by the multiplication of a function nonzero at $Q$ whose pullback is nonzero at $P$, so the valuation does not change. For any nonzero rational function $f$, write $f= s_Q^k\cdot u$. Then we have
\[
	\psi^*(f) = (\psi^*s_Q)^k\cdot (\psi^*u).
\]
The former has valuation $k\cdot e(P/Q)$ and the latter is nonzero at $P$. Hence we are done.
\epf

We also need the analytic interpretation of a local parameter. The
convergence in the following statement is used in
Section~\ref{sec:escaping-curves}.

\begin{lemma}[{\cite[Chapter 2, Sections 2.1--2.3]{shafarevich2013basic1} and \cite[Chapter 8, Section 1.1]{shafarevich2013basic2}}]
\label{lem: from algebra to analysis}
Let $\widetilde C$ be a smooth projective curve over $\mathbb C$, let
$P\in\widetilde C$, and let $s$ be a local parameter at $P$.
In the complex topology, there are a neighborhood $V$ of $P$ and
$\varepsilon>0$ such that
\[
 s|_V:V\longrightarrow D_\varepsilon
 :=\{z\in\mathbb C:|z|<\varepsilon\}
\]
is biholomorphic and sends $P$ to $0$. We may take $V$ inside a
Zariski-open neighborhood on which $s$ is regular. For every nonzero
$f\in\mathbb C(\widetilde C)$, after possibly shrinking this disc,
\[
 f\circ(s|_V)^{-1}(z)=\sum_{j=k}^{\infty}c_jz^j,
 \qquad k=\nu_P(f),\quad c_k\neq0,
\]
where the series converges on the punctured disc.
The amount of shrinking may depend on $f$.
\end{lemma}

Indeed, smooth complex algebraic curves are locally complex manifolds,
and algebraic local parameters give analytic coordinates, as explained in
\cite[Chapter 2, Section 2.3]{shafarevich2013basic1}.
Regular functions are holomorphic in these coordinates.
Writing $f=s^k h$ and expanding the holomorphic unit $h$ gives the
asserted Laurent expansion. For any fixed finite collection of rational
functions, one disc may be chosen for all of them by shrinking it once.

Let $P$ now be an ideal point of $C$ on its smooth projective model.
The birational map $\tilde\phi$ is regular and injective on a nonempty
Zariski-open subset of $\widetilde C$, whose complement is finite
\cite[Chapter I, Example 1.32 and Theorem 1.19]{shafarevich2013basic1}.
We can therefore choose $V$ above so that $\tilde\phi$ is regular and
injective on $V\setminus\{P\}$.

\begin{definition}\label{defn: puiseux branch}
	Suppose $C$ is an affine algebraic curve, and $(\widetilde{C},\tilde{\phi})$ is its smooth projective model. A \emph{Puiseux branch} is the choice of an ideal point $P\in \widetilde{C}$ and a local parameter $s$ at $P$, together with suitable $V$ and $\varepsilon$ supplied by the above lemma such that $V$ contains only one ideal point, namely $P$, and $\tilde{\phi}$ is injective on $V\backslash\{P\}$. Then we obtain the \emph{Puiseux parametrization}
	\[
		\tilde{\phi}\circ (s|_{V\backslash\{P\}})^{-1}: D_{\varepsilon}\backslash\{0\}\to C.
	\]
	Through $s|_V: V \to D_{\varepsilon}$, we may treat $s$ as the coordinate on $D_{\varepsilon}$. We treat the restriction
	\[
		\tilde{\phi}\circ (s|_{V\backslash\{P\}})^{-1}|_{D_{\varepsilon'}\backslash\{0\}}:D_{\varepsilon'}\backslash\{0\}\to C
	\]
	for $0<\varepsilon'<\varepsilon$ as the same Puiseux parametrization as the original.
\end{definition}

This parametrization is injective and holomorphic. In an affine
embedding of $C$, its coordinate functions have the convergent Laurent
expansions of Lemma~\ref{lem: from algebra to analysis}.
A proper closed subset of an irreducible curve is finite. An infinite
subset is therefore dense. If it is constructible, it
contains a nonempty open subset by the fact recalled in
Section~\ref{subsec:iteration}. We use this fact to pass from branches to
open subsets in Section~\ref{sec:escaping-curves}.

\subsection{Character varieties and peripheral trace functions}\label{sec: SL_2 character variety}
For a knot $K\subset S^3$, define
\[
	R_{\SL_2(\mathbb{C})}(K) = \operatorname{Hom}(\pi_K, \SL_2(\mathbb{C})).
\]
A finite presentation of $\pi_K$ realizes this as an affine algebraic set: the matrices assigned to the generators satisfy the determinant-one conditions and the equations given by the relators. Inversion is polynomial on $\SL_2(\mathbb{C})$.

\begin{definition}\label{defn: character and trace function}
The \emph{character} of $\rho\in R_{\SL_2(\mathbb{C})}(K)$ is the function $\chi_\rho(g)=\operatorname{tr}(\rho(g))$ on $\pi_K$. For $g\in \pi_K$, the associated \emph{trace function} is $t_g(\rho)=\operatorname{tr}(\rho(g))$ on the representation set. We denote by $T_{\SL_2(\mathbb{C})}(K)$ the subalgebra of $\mathbb{C}[R_{\SL_2(\mathbb{C})}(K)]$ generated by these functions.
\end{definition}

Note it is straightforward to check that all $t_g$ are regular functions. We now construct the character variety $X_{\SL_2(\mathbb{C})}(K)$ of $K$. Observe that the group $\SL_2(\mathbb{C})$ acts on $R_{\SL_2(\mathbb{C})}(K)$ by conjugation: for an element $h\in \SL_2(\mathbb{C})$ and $\rho \in R_{\SL_2(\mathbb{C})}(K)$, define
\[
	(h\cdot\rho)(g) = h\cdot \rho(g)\cdot h^{-1},~\forall g\in \pi_K
\]
Since $\SL_2(\mathbb{C})$ is reductive\footnote{We refer the readers to \cite{heusenerPorti2004PSL2Characters} for the definition of reductivity and the fact that $\SL_2(\mathbb{C})$ is reductive. This definition is different from the one in \cite[Definition 1.20]{brion2010actions}, which we used to show the GIT quotient map is surjective, but the equivalence follows from \cite[Definition 1.21 and Theorem 1.23]{brion2010actions}.}, it follows from \cite[Theorem 1.24]{brion2010actions}\footnote{The fact that the adjoint action is regular, which is needed as a hypothesis of the cited theorem, is straightforward to verify.} that $\mathbb{C}[R_{\SL_2(\mathbb{C})}(K)]^{\SL_2(\mathbb{C})}$ is finitely generated and defines an affine quotient with a surjective regular quotient map.

\begin{definition}\label{defn: character variety}
The $\SL_2(\mathbb{C})$ \emph{character variety} of $K$ is the affine quotient
\[
 X_{\SL_2(\mathbb{C})}(K)=R_{\SL_2(\mathbb{C})}(K)//\SL_2(\mathbb{C}).
\]
Its quotient map $\tau$ identifies its coordinate ring with the invariant ring:
\[
 \tau^*:\mathbb{C}[X_{\SL_2(\mathbb{C})}(K)]
 \xrightarrow{\cong}
 \mathbb{C}[R_{\SL_2(\mathbb{C})}(K)]^{\SL_2(\mathbb{C})}.
\]
\end{definition}

The term character variety does not require this affine algebraic set to be irreducible.

\begin{lemma}[{Sikora \cite[Proposition 2]{sikora2013generating}}]\label{lem: trace functions are all}
	Suppose $K\subset S^3$ is a knot. Then we have
	\[
		\mathbb{C}[X_{\SL_2(\mathbb{C})}(K)] \cong \mathbb{C}[R_{\SL_2(\mathbb{C})}(K)]^{\SL_2(\mathbb{C})} = T_{\SL_2(\mathbb{C})}(K)
	\]
	where $T_{\SL_2(\mathbb{C})}(K)$ is the trace algebra as in Definition~\ref{defn: character and trace function}.
\end{lemma}

\brem
A main ingredient of the proof of Proposition \ref{prop: 2bridge, valuation on I_p,q} below is the work by Ohtsuki \cite{ohtsuki1994ideal}. Ohtsuki's work then depends on a source paper by Culler--Shalen \cite{cullerShalen1983varieties}. In \cite{cullerShalen1983varieties}, the character variety $X_{\SL_2(\mathbb{C})}(K)$ is defined essentially as above, but with the trace algebra $T_{\SL_2(\mathbb{C})}(K)$, rather than the $\SL_2(\mathbb{C})$-invariant regular functions. We remark that Culler--Shalen's original construct coincides with what we defined in Definition~\ref{defn: character variety}, as a corollary to Lemma~\ref{lem: trace functions are all}.
\erem

\blem\label{lem: SL_2, character in bijection to points of character variety}
Suppose $K\subset S^3$ is a knot. Then there is a canonical bijection between the set $X_{\SL_2(\mathbb{C})}(K)$ and the set of all characters $\{\chi_{\rho}\}_{\rho\in R_{\SL_2(\mathbb{C})}(K)}$.
\elem
\bpf
Since $\tau$ is surjective, every point is $\tau(\rho)$ for some representation $\rho$. Two representations have the same image under $\tau$ precisely when all invariant regular functions agree on them. By Lemma~\ref{lem: trace functions are all}, this is equivalent to the equality of all their traces, that is, to equality of their characters.
\epf

With the above lemma at hand, we may write an element of the character variety $X_{\SL_2(\mathbb{C})}(K)$ as $\chi_{\rho}$. In particular, we have the following definition.
\begin{definition}\label{defn: I_p,q}
	Suppose $K\subset S^3$ is a knot. Use the fixed peripheral system $(m_K,\ell_K)$ with the preferred zero-linking longitude. We regard $m_K$ and $\ell_K$ as commuting elements in the knot group $\pi_K$. For $(p,q)\in\mathbb{Z}^2$, we define
	\[
		I_{p,q} = \left(t_{m_K^p\ell_K^q}\right)^2 - 4.
	\]
	Note the trace functions are regular functions on $X_{\SL_2(\mathbb{C})}(K)$ and so is $I_{p,q}$.
\end{definition}

\subsubsection{\texorpdfstring{$\PSL_2(\mathbb{C})$}{PSL2(C)} character variety}
By definition $\PSL_2(\mathbb{C}) = \SL_2(\mathbb{C})/\{\pm1\}$, and there is an algebraic group isomorphism $\PSL_2(\mathbb{C})\cong \SO_3(\mathbb{C})$. This isomorphism provides us with the facts we will need to construct the character variety: $\PSL_2(\mathbb{C})$ is reductive and $R_{\PSL_2(\mathbb{C})}(K)$ is indeed an affine algebraic set (cf. \cite[Section 2]{heusenerPorti2004PSL2Characters}).
\begin{definition}\label{defn: PSL_2C character variety}
	Suppose $K\subset S^3$ is a knot. Define
	\[
		R_{\PSL_2(\mathbb{C})}(K) = \operatorname{Hom}(\pi_K, \PSL_2(\mathbb{C})).
	\]
	Define its $\PSL_2(\mathbb{C})$ \emph{character variety} $X_{\PSL_2(\mathbb{C})}(K)$ to be the affine algebraic set arising from a construction similar to that in Section~\ref{sec: SL_2 character variety}, which has coordinate ring
	\[
		\mathbb{C}[X_{\PSL_2(\mathbb{C})}(K)] \cong \mathbb{C}[R_{\PSL_2(\mathbb{C})}(K)]^{\PSL_2(\mathbb{C})},
	\]
	and admits a surjective regular map $\hat{\tau}: R_{\PSL_2(\mathbb{C})}(K) \to X_{\PSL_2(\mathbb{C})}(K)$.
\end{definition}

\begin{definition}
	Suppose $K\subset S^3$ is a knot. For $\rho\in R_{\PSL_2(\mathbb{C})}(K)$ define the \emph{character} associated to it as
	\[
		\chi_{\rho}: \pi_K\to \mathbb{C}
	\]
	where $\chi_{\rho}(g) = \operatorname{tr}(\rho(g))^2$. Note $\operatorname{tr}(\rho(g))$ itself is not defined. More precisely, we write $\operatorname{tr}(\rho(g))^2 = \operatorname{tr}(A)^2$ where $A$ is a lift of $\rho(g)$ in $\SL_2(\mathbb{C})$. This lift is determined up to a sign but taking the square resolves the ambiguity. Similarly, for an element $g\in \pi_K$ we define the \emph{square trace function} associated to $g$ as $\hat{t}_g(\rho) = \operatorname{tr}(\rho(g))^2$. Clearly, $\hat{t}_g$ are $\PSL_2(\mathbb{C})$ invariant and hence become regular functions on $X_{\PSL_2(\mathbb{C})}(K)$.
\end{definition}

\begin{lemma}[{Heusener--Porti \cite[Theorem 1.3 and Proposition 2.2]{heusenerPorti2004PSL2Characters}}]\label{lem: square trace function generates coordinate ring}
	The square trace functions $\hat{t}_g$ for all $g\in \pi_K$ generate $\mathbb{C}[R_{\PSL_2(\mathbb{C})}(K)]^{\PSL_2(\mathbb{C})}$. Consequently, there is a (canonical) bijection between the set of all characters and the set $X_{\PSL_2(\mathbb{C})}(K)$, which identifies $\chi_{\rho}$ with $\hat{\tau}(\rho)$ for any $\rho\in R_{\PSL_2(\mathbb{C})}(K)$.
\end{lemma}

For the comparison below, we use finite morphisms between affine algebraic sets, without requiring these sets to be irreducible. We need to align the terminology in this paper with that in \cite{ferrersantos2005actions1}: (1) Affine algebraic sets as defined in Section~\ref{subsec:iteration} are algebraic varieties in \cite[Chapter 1, Observations 4.19 and 4.39]{ferrersantos2005actions1}. (2) A regular map between two affine algebraic sets as in Section~\ref{subsec:iteration} is a morphism between algebraic varieties as in \cite[Chapter 1, Definitions 3.28 and 4.37, Observation 3.29, and Corollary 3.38]{ferrersantos2005actions1}. (3) Such a map is finite in the sense of \cite[Definition 4.74]{ferrersantos2005actions1} if the coordinate ring of its source is finitely generated as a module over the coordinate ring of its target via pullback; the global coordinate rings then provide the required trivial affine cover.

\begin{lemma}\label{lem: surjective finite morphism from SL_2 to PSL_2}
	Suppose $K\subset S^3$ is a knot. Then there is a surjective regular map (cf. Section~\ref{subsec:iteration}) $\psi: X_{\SL_2(\mathbb{C})}(K) \to X_{\PSL_2(\mathbb{C})}(K)$ such that the coordinate ring $\mathbb{C}[X_{\SL_2(\mathbb{C})}(K)]$ is a finitely generated $\mathbb{C}[X_{\PSL_2(\mathbb{C})}(K)]$-module via the pullback
	\[
		\psi^*: \mathbb{C}[X_{\PSL_2(\mathbb{C})}(K)]\to \mathbb{C}[X_{\SL_2(\mathbb{C})}(K)].
	\]
\end{lemma}
\bpf
First note that for the knot group $\pi_K$, we have $H^2(\pi_K;\mathbb{Z}/2\mathbb{Z}) = 0$ (cf. the discussion before \cite[Lemma 3.3]{boyerZhang1998CullerShalen}), and hence any $\PSL_2(\mathbb{C})$ representation lifts to an $\SL_2(\mathbb{C})$ representation (cf. \cite[Section 4]{heusenerPorti2004PSL2Characters}). Then Heusener--Porti \cite[Proposition 4.2]{heusenerPorti2004PSL2Characters} supply the identification of $X_{\PSL_2(\mathbb{C})}(K)$ with the GIT quotient
\[
	X_{\PSL_2(\mathbb{C})}(K) \cong X_{\SL_2(\mathbb{C})}(K) // H^1(\pi_K;\mathbb{Z}/2\mathbb{Z}).
\]
Note that $H^1(\pi_K;\mathbb{Z}/2\mathbb{Z}) = \operatorname{Hom}(\pi_K,\mathbb{Z}/2\mathbb{Z}) = \mathbb{Z}/2\mathbb{Z} = \{\pm1\}$. So we obtain a regular map
\[
	\psi: X_{\SL_2(\mathbb{C})}(K) \to X_{\PSL_2(\mathbb{C})}(K)
\]
whose pullback provides an identification 
\[
	\psi^*: \mathbb{C}[X_{\PSL_2(\mathbb{C})}(K)] \cong \mathbb{C}[X_{\SL_2(\mathbb{C})}(K)]^{\{\pm1\}}
\] 
Then $\mathbb{C}[X_{\SL_2(\mathbb{C})}(K)]$ is indeed a finitely generated module over $\mathbb{C}[X_{\SL_2(\mathbb{C})}(K)]^{\{\pm1\}}$: let $t_{g_1},\cdots,t_{g_n}$ be a set of ($\SL_2(\mathbb{C})$) trace functions that generate $\mathbb{C}[X_{\SL_2(\mathbb{C})}(K)]$ as a $\mathbb{C}$-algebra. Then we know $t_{g_i}^2\in \mathbb{C}[X_{\SL_2(\mathbb{C})}(K)]^{\{\pm 1\}}$, and hence $t_{g_1}^{\varepsilon_1}\cdots t_{g_n}^{\varepsilon_n}$ for all $\varepsilon_i\in\{0,1\}$ generate $\mathbb{C}[X_{\SL_2(\mathbb{C})}(K)]$ as a $\mathbb{C}[X_{\SL_2(\mathbb{C})}(K)]^{\{\pm1\}}$-module.
\epf

\begin{remark}
	We now provide a concrete description of the map $\psi$ constructed above. From Lemma~\ref{lem: SL_2, character in bijection to points of character variety} and Lemma~\ref{lem: square trace function generates coordinate ring}, we can identify points in the character variety with characters. If $\rho\in R_{\SL_2(\mathbb{C})}(K)$ is an $\SL_2(\mathbb{C})$ representation, it naturally induces a $\PSL_2(\mathbb{C})$ representation $\hat{\rho}$ by composing with the natural quotient map $\SL_2(\mathbb{C})\to \PSL_2(\mathbb{C})$. Then in terms of characters, $\psi(\chi_{\rho}) = \chi_{\hat{\rho}}$.
\end{remark}

\begin{corollary}\label{cor: C to D}
	Suppose $D\subset X_{\PSL_2(\mathbb{C})}(K)$ is a one-dimensional irreducible Zariski-closed subset. Then there exists a one-dimensional irreducible Zariski-closed subset $C\subset X_{\SL_2(\mathbb{C})}(K)$ such that
	\[
		\psi|_C : C \to D
	\]
	is a surjective morphism.
\end{corollary}
\bpf
Let $C_1,\cdots,C_n$ be the irreducible components of $\psi^{-1}(D)$ (by \cite[Chapter I, Corollary 1.6]{hartshorne1977algebraic} there are finitely many). By Lemma~\ref{lem: surjective finite morphism from SL_2 to PSL_2} and the convention above, $\psi$ is finite.

Thus, we apply \cite[Chapter 1, Theorem 4.93]{ferrersantos2005actions1} and conclude that $\psi$ is a closed map, so each $\psi(C_i)$ is a closed subset of $D$. By the definition of irreducibility in Section~\ref{subsec:iteration}, there exists an index $i$ such that $\psi(C_i) = D$. Hence take $C = C_i$. It remains to show that $C$ is indeed a curve. By \cite[Chapter 1, Theorem 4.93]{ferrersantos2005actions1}, $\psi|_C$ is a finite map. Then \cite[Chapter 1, Lemma 4.95]{ferrersantos2005actions1} implies that $C$ must be a curve.
\epf

\begin{remark}
    $I_{p,q}$ was originally defined on $X_{\SL_2(\mathbb{C})}(K)$ as an element in the coordinate ring. However, Lemma~\ref{lem: surjective finite morphism from SL_2 to PSL_2} canonically identifies $\mathbb{C}[X_{\PSL_2(\mathbb{C})}(K)]$ as the invariant subring $\mathbb{C}[X_{\SL_2(\mathbb{C})}(K)]^{\{\pm1\}}$, and $I_{p,q}$ lies inside that subring. So, abusing notation, we write $I_{p,q}$ as a polynomial function on $X_{\PSL_2(\mathbb{C})}(K)$ as well, and we have $\psi^*I_{p,q} = I_{p,q}$ where the function $I_{p,q}$ on the left-hand side is a polynomial function on $X_{\PSL_2(\mathbb{C})}(K)$ while the one on the right-hand side is on $X_{\SL_2(\mathbb{C})}(K)$.
\end{remark}

\subsection{Ohtsuki's result}\label{sec: ohtsuki's result}
In this section, we summarize key ingredients from Ohtsuki \cite{ohtsuki1994ideal}. Let $K$ be a $2$-bridge knot, then its fundamental group has a presentation of the form
\[
	\pi_K = \langle u,v~|~wu=vw\rangle
\]
for a suitable word $w$ in $u$ and $v$. Ohtsuki defines $R$ to be the dimension $\geq 4$ components of the variety of all $\PSL_2(\mathbb{C})$ representations, and $X$ to be the set of all conjugacy classes of elements in $R$. For a generic element $[\rho]\in X$, he defined 
\[
	\tau = \frac{2\operatorname{tr}(\widetilde U\widetilde V)-(\operatorname{tr}\widetilde U)^2}{(\operatorname{tr}\widetilde U)^2-4} \text{ and } \xi = \frac{-2}{({\rm tr}(\widetilde{U}))^2 - 4}.
\] 
Here we take $U = \rho(u)$, $V = \rho(v)$ and $\SL_2(\mathbb{C})$-lifts $\widetilde{U}$ and $\widetilde{V}$ which are conjugate to each other. Ohtsuki showed that these coordinates identify the set of generic irreducible points (i.e., conjugacy classes) of $X$ with a Zariski-open subset $\mathcal{C}(K)$ of the one-dimensional algebraic set in $\mathbb C^2$ cut out by $z(\tau,\xi)=0$. The following lemma is implicit in his paper.
\begin{lemma}\label{lem: square trace are regular in Ohsuki's setup}
	For any element $g\in \pi_K$, the map $\hat{t}_g:\mathcal{C}(K)\to \mathbb{C}$ defined by
	\[
		\hat{t}_g([\rho]) = ({\rm tr}(\rho(g)))^2
	\]
	is regular on $\mathcal{C}(K)$.
\end{lemma}
The above lemma is not explicitly proved, but is essentially a corollary to \cite[Lemma 1.2]{ohtsuki1994ideal} with proof similar to that of \cite[Lemma 4.9]{ohtsuki1994ideal}. As a corollary, Ohtsuki stated that functions
\[
	I^{Oh}_{p,q}([\rho]) = \bigg({\rm tr}(\rho(m_K^p\ell_K^q))\bigg)^2 - 4.
\]
for all $(p,q)\in\mathbb{Z}^2$ are all regular on $\mathcal{C}(K)$. The key ingredient is the following.
\begin{lemma}\label{lem: ohtsuki's pole order}
Suppose $K\subset S^3$ is a $2$-bridge knot and $r\neq 0$ is the boundary slope of a non-closed oriented essential (i.e., incompressible, boundary-incompressible, and non-boundary-parallel) properly embedded orientable surface in the knot complement, then there exists an irreducible component ${X}_0$ of $\mathcal{C}(K)$, a smooth projective model $\widetilde{X}_0$ with a point $Q\in\widetilde{X}_0$, and $\varepsilon_0\in\{1,2\}$ such that for any $(p,q)\in \mathbb{Z}^2$, we have
	$${\rm pole}_Q(I^{Oh}_{p,q}|_{{X}_0}) = \varepsilon_0\cdot |p - rq|.$$
\end{lemma}

\begin{proof}
	By Hatcher–Thurston’s classification (cf. \cite[Theorem 1]{hatcherThurston1985incompressible}), as recalled in~\cite[Section 2]{ohtsuki1994ideal}, there is a continued fraction $[n_j]$, with $|n_j|\ge2$ for every $j$, such that $r = N_{[n_j]}$ (cf. \cite[Section 2]{ohtsuki1994ideal} for definition). The unique continued fraction expansion $[n_j]$ whose coefficients are all even has corresponding $N_{[n_j]} = 0$. Since $r\ne0$, in an expansion $[n_j]$ such that $r = N_{[n_j]}$, some $n_j$ is odd. Consequently, choosing each residue $k_j\in\mathbb Z/(n_j)$ to be the class of $1$ gives a tuple allowed by Ohtsuki’s Theorem 5.1. Thus we can apply that theorem and the pole-order formula immediately preceding Proposition 5.2 and take $X_0$, $Q$ and $\varepsilon_0\in\{1,2\}$ as desired.
\end{proof}

\section{Two-bridge knots and their Whitehead doubles}
\subsection{Whitehead-iterable branches}\label{sec:escaping-curves}

For every finite $N$, the set $\mathcal S_N$ of
Section~\ref{subsec:iteration} contains a dense open subset of the
peripheral torus. We now give a condition ensuring that
a curve of actual augmentation points meets these sets. The condition is
expressed in terms of valuations on a smooth projective model, using the
terminology of the preceding section.

We work over $\mathbb C$ and write
\[
 T=(\mathbb C^*)^2_{(\mu,\eta)},\qquad \eta= \lambda^2.
\]
We use $\mathcal S_N$ and $\mathcal G_N$ for the sets already defined in Section \ref{subsec:iteration}; in particular, $\mathcal G_N\subseteq\mathcal S_N$.

\begin{definition}\label{def:escaping-branch}
Let $C\subset T$ be an irreducible algebraic curve. We say that $C$ has a
\emph{Whitehead-iterable branch} (and iterable branch for short), if $C$ admits an ideal point $P$ with $\nu_P(\mu) = -p$ and $\nu_P(\eta) = q$ for $(p,q)\in\mathbb{Z}^2$ such that
\begin{equation}\label{eq: iterable condition}
	p>0 \text{ and (either }q>p~{\rm or}~q<0).
\end{equation}
\end{definition}

Now for each point $(\mu,\eta)$, we have one-step correspondence $(\widetilde{\mu},\widetilde{\lambda})$ as in \eqref{eq:raw-whitehead-correspondence}. Since $\widetilde{\mu}$ is defined via a quadratic equation, in general it is not a rational function on the original $C$, so the valuation $\nu$ is not directly applicable for $\widetilde{\mu}$. There are two routes to settle this issue: (1) Take a suitable (branched) double cover of $C$ so that $(\widetilde{\mu},\widetilde{\lambda})$ are well defined rational functions, and work with valuations of them, and (2) switch to use complex analytical topology. The two routes will lead to equivalent arguments, and we adopt (2). 

As in Definition \ref{defn: puiseux branch}, in a small enough complex-analytical neighborhood of the ideal point $P$ of $C$, we have a local Puiseux parameter $s$. Thus we introduce the following notation for a local (near $P$) meromorphic function $f$ with at most a pole at $s=0$: one can express such an $f$ that is not identically zero as
\[
	f = s^r\cdot h(s)
\]
where $r\in\mathbb{Z}$ and $h$ is holomorphic with $h(0)\neq 0$. Then we write $\mathfrak{o}(f) = r$. Also define $\mathfrak{o}(0) = +\infty$. Note $\mu$ and $\eta$ restrict to meromorphic functions near $s = 0$, and the above valuation, via Lemma \ref{lem: from algebra to analysis}, can be restated as 
\[
	\mathfrak{o}(\mu) = -p \aand \mathfrak{o}(\eta) = q.
\]
Furthermore, we have the following properties for $\mathfrak{o}$:
\begin{itemize}
	\item $\mathfrak{o}(c) = 0$ for $c\in\mathbb{C}^*$.
	\item $\mathfrak{o}(f\cdot g) = \mathfrak{o}(f) + \mathfrak{o}(g)$.
	\item $\mathfrak{o}(f+g) \geq \min\{\mathfrak{o}(f),\mathfrak{o}(g)\}$, with equality whenever $\mathfrak{o}(f)\neq \mathfrak{o}(g)$.
\end{itemize}

\begin{theorem}\label{thm:escaping-nontrapping}
Let $C\subset T$ be an irreducible algebraic curve with an iterable branch. For every integer $N\geq0$, the intersection $C\cap\mathcal S_N$ (cf. Section \ref{subsec:iteration})
contains a nonempty Zariski-open subset of $C$.
\end{theorem}

\begin{proof}
Let $P\in\widetilde C$ be the ideal point specified by the
iterable branch, and let $u$ be a local parameter at $P$. On a sufficiently small parameter, write $\eta(u)=u^q h(u)^2$, where $h$ is holomorphic near $0$
and $h(0)\neq0$. On the covering disk given by $u=s^2$, define
\[
 \mu_0(s):=\mu(s^2),\qquad
 \lambda_0(s):=s^q h(s^2),\qquad
 \eta_0(s):=\lambda_0(s)^2=\eta(s^2).
\]
Throughout this proof, $\mathfrak{o}$ denotes the order at $s=0$ on this covering disk. Thus
\[
 -\mathfrak{o}(\mu_0)=2p,\qquad
 \mathfrak{o}(\eta_0)=2q.
\]
We will establish the following claim by induction via the correspondence $\mathcal W$ introduced in Section \ref{subsec:iteration}.

\noindent{\bf Claim}. For any $N\in\mathbb{Z}_{\geq 0}$, there exists $\varepsilon_N>0$ such that on $D_{\varepsilon_N} = \{z\in\mathbb{C}~|~|z|<\varepsilon_N\}$,    there exist meromorphic functions $\mu_N(s)$, $\lambda_N(s)$ with only possible pole at $s=0$,
\[
	\eta_N := \lambda_N^2,~ p_N := -\mathfrak{o}(\mu_N),\aand q_N := \mathfrak{o}(\eta_N)
\]
such that the following conditions hold.
\begin{itemize}
	\item We have $\mu_0 = \mu(s^2)$ and $\eta_0 = \eta(s^2)$.
	\item We have the following inductive relation for any $N\in\mathbb{Z}_{\geq 0}$:
	\[
		p_{N+1} = p_{N}\aand q_{N+1} = \begin{cases}
		2(p_{N}+q_{N}) & q_{N} > p_{N}\\
		2q_{N} & q_{N} < 0
	\end{cases}
	\]
	As a consequence, all $(p_N,q_N)$ satisfy \eqref{eq: iterable condition}.
	\item We have $(\mu_{N}(s),\eta_{N}(s),\mu_{N+1}(s),\lambda_{N+1}(s))\in \mathcal{W}_{\rm step}$ (cf. \eqref{eq:good-whitehead-correspondence}) for all $0<|s|<\min\{\varepsilon_N, \varepsilon_{N+1}\}$. As a consequence, for all $N\in\mathbb{Z}_{\geq 0}$, $(\mu_0(s),\eta_0(s))\in C\cap\mathcal G_N$ for all $0<|s|<\min\{\varepsilon_0,\cdots,\varepsilon_N\}$.
\end{itemize}

We first note that the above claim suffices to prove the theorem. Indeed, this claim guarantees that there are infinitely many points
of $C\cap\mathcal G_N$. By
Proposition~\ref{prop:generic-finite-depth-chains}, this intersection is
constructible. An infinite constructible subset of an irreducible
algebraic curve contains a nonempty Zariski-open subset. Finally, the fact that
$\mathcal G_N\subseteq\mathcal S_N$ as in \eqref{eq: G inside S} gives the assertion.

We now proceed to prove the claim. Induction starts with $(\mu_0(s),\eta_0(s))$. Since $p_0=2p$ and $q_0=2q$, the pair $(p_0,q_0)$ satisfies \eqref{eq: iterable condition}. Assume $\mu_n(s),\eta_n(s)$ have been constructed with the desired properties in the claim. We now build $\mu_{n+1}$ and $\eta_{n+1}$.
Pointwise, $\mu_{n+1}(s)$ is a solution to the equation
\[
	A_n(s)z^2 + B_n(s)z + A_n(s) = 0.
\]
where $A_n$ and $B_n$ are derived from \eqref{eq:determinant-polynomial} by inserting $\mu_n(s)$ and $\eta_n(s)$. Explicitly, we have
\[
 \begin{aligned}
 A_n&=-\mu_n(\eta_n-1)(\mu_n\eta_n+1),\\
 B_n&=\mu_n^3\eta_n^2+2\mu_n^2\eta_n-\mu_n^2+\mu_n\eta_n^2-2\mu_n\eta_n-1.
 \end{aligned}
\]

They are obviously meromorphic on a sufficiently small disk. A direct computation of the order yields 

\begin{table}[http]
%\color{orange}
\centering
\begin{tabular}{|c|c|c|c||c|c|c|c|c|}
\hline
&$\mu_n$&$\eta_n-1$&$\mu_n\eta_n + 1$&$\mu_n^3\eta_n^2$&$\mu_n^2\eta_n$&$\mu_n^2$&$\mu_n\eta_n^2$&$\mu_n\eta_n$ \\
\hline
$q_n>p_n$&$-p_n$&0&0&$2q_n-3p_n$&$q_n-2p_n$&$\mathbf{-2p_n}$&$2q_n-p_n$&$q_n-p_n$ \\
\hline
$q_n<0$&$-p_n$&$q_n$&$q_n-p_n$&$\mathbf{2q_n-3p_n}$&$q_n-2p_n$&$-2p_n$&$2q_n-p_n$&$q_n-p_n$ \\
\hline
\end{tabular}
\par\vspace{6pt}
\caption{$\mathfrak{o}$ of building blocks of $A_n$ and $B_n$. ``Leading'' terms for $B_n$ are bolded.}
\label{tab: o of A and B}
\end{table}
As a result, we have
\[
\mathfrak{o}\left(\frac{A_n}{B_n}\right) = \begin{cases}
	(-p_n) - (-2p_n) = p_n & q_n > p_n\\
	(2q_n-2p_n) - (2q_n-3p_n) = p_n & q_n < 0
\end{cases}
\]
As a conclusion, we obtain an equation
\[
\frac{A_n}{B_n}z^2 + z + \frac{A_n}{B_n} = 0
\]
with holomorphic coefficients. At $s = 0$, we must have $z=0$ and derivative of the left-hand-side at $z = 0$ is $1\neq 0$. Thus the holomorphic implicit function theorem applies and there exists a local (non-constant) holomorphic function $z_{n+1}(s)$ satisfying the equation. We then take $\mu_{n+1}(s) = z_{n+1}^{-1}(s)$, as a meromorphic function for sufficiently small $|s|$. From Vieta's theorem, we know
\[
	z_{n+1}(s) + \mu_{n+1}(s) = -\frac{B_n(s)}{A_n(s)}
\]
Hence we conclude (note we must have $\mathfrak{o}(\mu_{n+1}) < 0 < \mathfrak{o}(z_{n+1})$)
\[	
	\mathfrak{o}(\mu_{n+1}) = \mathfrak{o}(z_{n+1} + \mu_{n+1}) = \mathfrak{o}\left(\frac{B_n}{A_n}\right) = -p_n.
\]
$\eta_{n+1}=\lambda_{n+1}^2$ is uniquely determined by $\mu_n$, $\eta_n$, and $\mu_{n+1}$ via \eqref{eq:raw-whitehead-correspondence}. However, we need a different formula to perform the computation.

Recall, as in the proof of Lemma \ref{lem:whitehead-block-existence}, we have
\[
	d_n=\frac{1+\mu_{n+1}}{1+\mu_n}\aand a_n=\mu_n d_n.
\]
Definition \ref{def:admissible-doubling-chain} item (2) implies an identity
\begin{equation}\label{eq: d_n - 1}
	(\mu_n\eta_n^2-1)(a_n-1)(d_n-1)=\eta_n(\mu_n-1)a_nd_n.
\end{equation}
Direct computation shows
\[
\mathfrak{o}(d_n) = 0,~ \mathfrak{o}(a_n) = -p_n,~\mathfrak{o}(\mu_n-1) = -p_n, ~\mathfrak{o}(a_n-1) = -p_n,
\]
and
\[
\mathfrak{o}(\mu_n\eta_n^2 - 1) =\begin{cases}
	0 & q_n > p_n \\
	2q_n - p_n & q_n < 0. 
\end{cases}
\]
Hence we conclude
\[
\mathfrak{o}(d_n - 1) = \begin{cases}
 q_n - p_n & q_n > p_n \\
 -q_n & q_n < 0.
\end{cases}
\]
Definition \ref{def:admissible-doubling-chain} item (3) implies
\begin{equation}\label{eq: expression for eta_n+1}
	\eta_{n+1} = \left(\frac{\eta_n(a_n-1)+(d_n-1)}{(a_n-1)+\eta_n(d_n-1)}\right)^2.
\end{equation}
A direct computation for the denominator shows that $\mathfrak{o}((a_n-1)+\eta_n(d_n-1)) = -p_n$ for both cases. We now focus on the numerator. For the case $q_n < 0$, we compute $\mathfrak{o}(\eta_n(a_n-1)+(d_n-1)) = q_n - p_n$. For the case $q_n > p_n$, Equation \eqref{eq: d_n - 1} further implies
\[
	\begin{aligned}
		\eta_n(a_n-1)+(d_n-1) &= \eta_n(a_n-1) + \frac{\eta_n(\mu_n-1)a_nd_n}{(\mu_n\eta_n^2-1)(a_n-1)}\\
		&= \eta_n\frac{(\mu_n\eta_n^2-1)(a_n-1)^2+a_n\mu_n d_n - a_n d_n}{(\mu_n\eta_n^2-1)(a_n-1)}\\
        &= \eta_n\frac{\mu_n\eta_n^2(a_n-1)^2 - a_n^2 + 2a_n - 1 + a_n\mu_n d_n - a_n d_n}{(\mu_n\eta_n^2-1)(a_n-1)}\\
		\text{(recall }a_n=\mu_nd_n\text{) }&= \eta_n\frac{\mu_n\eta_n^2(a_n-1)^2+a_n(1+1-d_n)-1}{(\mu_n\eta_n^2-1)(a_n-1)}
	\end{aligned}
\]
Hence we conclude $\mathfrak{o}(\eta_n(a_n-1)+(d_n-1)) = q_n$. In summary, we conclude
\[
q_{n+1} = \mathfrak{o}(\eta_{n+1}) = \begin{cases}
	2(q_n + p_n) & q_n > p_n\\
	2q_n & q_n < 0
\end{cases}
\]
Finally, it remains to show that $(\mu_n(s),\eta_n(s),\mu_{n+1}(s),\lambda_{n+1}(s))\in \mathcal W_{\rm step}$ when $|s|$ is small enough and non-zero. This amounts to checking several non-vanishing conditions as in \eqref{eq:peripheral-quadratic-discriminant} and \eqref{eq:good-whitehead-correspondence}. The expression $\delta(\mu_n,\eta_n,\mu_{n+1})$ can be expressed as the denominator of $\lambda_{n+1}$ in \eqref{eq: expression for eta_n+1} times $-(\mu_n+1)$, which has order $-2p_n$ and hence non-zero for small enough $|s|$. Remaining terms are:
\[
\eta_n-1\neq 0,~\mu_n\pm1\neq 0,~\mu_n\eta_n+1\neq 0,~\mu_n^2\eta_n^2-\mu_n\eta_n^2+4\mu_n\eta_n-\mu_n+1\neq 0,~\mu_n\eta_n^2-1\neq 0.
\]
This can be verified by direct computation:
\begin{table}[http]
%\color{orange}
\centering
\begin{tabular}{|c|c|c|c|c|c|}
\hline
&$\eta_n - 1$&$\mu_n\pm1$&$\mu_n\eta_n + 1$&$\mu_n^2\eta_n^2-\mu_n\eta_n^2+4\mu_n\eta_n-\mu_n+1$&$\mu_n\eta_n^2-1$ \\
\hline
$q_n>p_n$&$0$&$-p_n$&0&$-p_n$&$0$ \\
\hline
$q_n<0$&$q_n$&$-p_n$&$q_n-p_n$&${2q_n-2p_n}$&$2q_n-p_n$ \\
\hline
\end{tabular}
\par\vspace{6pt}
\caption{$\mathfrak{o}$ of the five terms in the non-vanishing conditions.}
\label{tab: non-vanishing order}
\end{table}

In summary, each has finite order in both cases, and hence must be non-zero when $|s|$ is positive but small enough. We thus complete the induction and conclude the proof of the claim and hence the theorem.
\end{proof}

\begin{remark}
The passage to the covering disk $u=s^2$ ensures the existence of the initial meromorphic $\lambda_0$. Although $\eta$ is meromorphic in the original parameter $u$, its order $q$ may be odd, whereas the square of a meromorphic function has even order. After the substitution $u=s^2$, the function $\eta_0(s)=\eta(s^2)$ has order $2q$ and hence admits a meromorphic square root on a sufficiently small disk. If $q$ is even, this covering substitution can be omitted. This applies, for example, to the branch constructed for a $2$-bridge knot in
Theorem~\ref{thm: 2-bridge knot admits iterable branch}, since $q=rp$, where $r$ is an even integer.
\end{remark}

\begin{definition}\label{def:actual-kch-curve}
Let $K$ be a knot. An irreducible curve $C\subset T$ is an
\emph{actual rank-$d$ augmentation peripheral curve} for $K$ if it contains
a Zariski-dense subset $C^{\mathrm{act}}$ such that every point of
$C^{\mathrm{act}}$ is of the form
\[
 (\mu,\eta)=\bigl(\epsilon(\mu),\epsilon(\lambda)^2\bigr)
\]
for a rank-$d$ augmentation $\epsilon$ of $K$ with
$\epsilon(\mu)\neq1$.
\end{definition}

Under the augmentation–sheaf correspondence, these points are
$(\mu,\lambda^2)$, where $(\mu,\lambda)$ are the peripheral data of
irreducible $d$-dimensional diagonalizable $\murk$-one representations.
The word ``actual'' distinguishes these points from points introduced
only by taking a Zariski closure.

\begin{corollary}\label{cor:escaping-rank-growth}
Suppose that $K$ has an actual rank-$d$ augmentation peripheral curve with
an iterable branch. Then, for every $N\geq0$,
\[
 \ar(\Wh^N(K))\geq2^Nd.
\]
If $d=b(K)$, then
\[
 \ar(\Wh^N(K))=\mr(\Wh^N(K))=b(\Wh^N(K))=2^Nb(K).
\]
\end{corollary}

\begin{proof}
By Theorem~\ref{thm:escaping-nontrapping}, $C\cap\mathcal S_N$ contains
a nonempty open subset of $C$. This subset meets $C^{\mathrm{act}}$.
Choose a rank-$d$ augmentation giving such a point and apply
Proposition~\ref{thm:finite-depth-iteration}. This proves the lower bound.
When $d=b(K)$, Corollary~\ref{cor:finite-depth-rank-equality} gives all
three equalities.
\end{proof}

The actual augmentation point chosen in this proof may depend on $N$.
The result does not require a single point to realize all the iterations.

\subsection{Two-bridge knots admit iterable branches}\label{sec:slope-detection}
Throughout this section, we take $K$ to be a $2$-bridge knot. The goal of this section is to establish the following theorem
\begin{theorem}\label{thm: 2-bridge knot admits iterable branch}
	Any $2$-bridge knot admits an actual rank-$2$ augmentation peripheral curve as in Definition \ref{def:actual-kch-curve} with an iterable branch as in Definition \ref{def:escaping-branch}.
\end{theorem}

\begin{proof}[Proof of Theorem~\ref{thm:main-two-bridge}]
The combination of Theorem \ref{thm: 2-bridge knot admits iterable branch} and Corollary \ref{cor:escaping-rank-growth} directly implies
\[
 \ar(\Wh^N(K))=\mr(\Wh^N(K))=b(\Wh^N(K))=2^{N+1}
\]
for every $N\geq0$.
\end{proof}

The rest of the section is dedicated to the proof of Theorem \ref{thm: 2-bridge knot admits iterable branch}. We need several preparations. First, a properly embedded surface $S\subset E(K) = S^3 - {\rm int}(N(K))$ is said to be \emph{essential} if it is incompressible, boundary-incompressible, and non-boundary-parallel. We have the following existence result.
\begin{lemma}\label{lem: non-zero even boundary slope}
	Suppose $K\subset S^3$ is a $2$-bridge knot. Then there exists an oriented essential non-closed properly embedded surface $S$ in its knot complement and $\partial S$ has slope $r$ which is an even non-zero integer.
\end{lemma}
\begin{proof}
	First, \cite[Theorem 1]{mattmanMaybrunRobinson2008twoBridge} shows that the diameter of the set of boundary slopes is positive for a $2$-bridge knot $K$. Since a minimal genus Seifert surface has boundary slope $0$, there must exist a non-zero boundary slope $r$. Furthermore, \cite[Section 1]{mattmanMaybrunRobinson2008twoBridge} indicates that such an $r$ can be picked to be an even integer.
\end{proof}

\begin{proposition}\label{prop: 2bridge, valuation on I_p,q}
	Suppose $K\subset S^3$ is a nontrivial $2$-bridge knot. Then there exist a nonzero even boundary slope $r$, a one-dimensional irreducible Zariski-closed subset $C'\subset X_{\SL_2(\mathbb{C})}(K)$, an ideal point $Q$ of the smooth projective model $\widetilde{C}'$ of $C'$, and a positive rational number $\varepsilon$, such that for any $(p,q)\in \mathbb{Z}^2$ we have
\[
\operatorname{pole}_Q(I_{p,q}|_{C'}) = \varepsilon\cdot |p - rq|.
\]
Here, $I_{p,q}$ is a regular function on $X_{\SL_2(\mathbb{C})}(K)$ introduced in Definition~\ref{defn: I_p,q}, and we view $I_{p,q}|_{C'}$ as a rational function on $\widetilde{C}'$.
\end{proposition}
\begin{proof}
We first show that there exists a 1-dimensional irreducible Zariski-closed subset $D\subset X_{\PSL_2(\mathbb{C})}(K)$, and an ideal point $Q$ of its smooth projective model $\widetilde{D}$, such that there exists a positive integer $\varepsilon'$ and for any $(p,q)\in \mathbb{Z}^2$, we have
\[
	{\rm pole}_Q(I_{p,q}|_D) = \varepsilon'\cdot |p - rq|.
\]

We proceed with the non-zero even boundary slope $r$ obtained from Lemma \ref{lem: non-zero even boundary slope}. Lemma \ref{lem: ohtsuki's pole order} then applies and produces a quasi-affine curve $X_0$ with a point $Q$ in its smooth projective model $\widetilde{X}_0$ and the pole order formula
\begin{equation}\label{eq: ohtsuki's pole order formula}
	{\rm pole}_{Q}(I^{Oh}_{p,q}|_{X_0}) = \varepsilon_0\cdot |p - rq|.
\end{equation}

However, we still need a transition from Ohtsuki's space $X$ to our $\PSL_2(\mathbb{C})$ character variety. Ohtsuki defined elements in $X$ as conjugation classes of representations $\rho\in R_{\PSL_2(\mathbb{C})}(K)$, and our map $\hat{\tau} : R_{\PSL_2(\mathbb{C})}(K)\to X_{\PSL_2(\mathbb{C})}(K)$ defined as in Definition \ref{defn: PSL_2C character variety} with a chosen set of square trace functions generating the coordinate ring naturally induces a map, which, by abusing notation, we still denote by $\hat{\tau}: X_0 \to X_{\PSL_2(\mathbb{C})}(K)$. One core identity that is straightforward from the construction is that
\[
	I_{p,q}^{Oh} = I_{p,q}\circ \widehat{\tau}
\]

We claim that the map $\hat{\tau}$ is regular. Indeed, by Lemma \ref{lem: square trace function generates coordinate ring}, the coordinate ring $\mathbb{C}[X_{\PSL_2(\mathbb{C})}(K)]$ is generated by square trace functions. Lemma \ref{lem: square trace are regular in Ohsuki's setup} implies that all square trace functions are regular on $X_0$.

Next, take $D = \overline{\hat{\tau}(X_0)}$ to be the Zariski-closure of the image inside $X_{\PSL_2(\mathbb{C})}(K)$. It is irreducible since $X_0$ is. \cite[Corollary 3.6]{heusenerPorti2004PSL2Characters} states that the irreducible locus $X^{irr}_{\PSL_2(\mathbb{C})}(K) \subset X_{\PSL_2(\mathbb{C})}(K)$ is Zariski open. Taking $D^{irr} = X^{irr}_{\PSL_2(\mathbb{C})}(K)\cap D$, it is open in $D$. Thus, $\hat{\tau}^{-1}(D^{irr})$ is open, and we show that it is not empty in $X_0$. Suppose, on the contrary, that every element $[\rho]\in X_0$ is reducible. Then for any $[\rho]\in X_0$, $\hat{\tau}([\rho]) = \chi_{\rho}$ is a reducible character. \cite[Remark 3.7]{heusenerPorti2004PSL2Characters} implies that the same $\chi_{\rho}$ can also be represented by an abelian representation $\rho'$: $\chi_{\rho} = \chi_{\rho'}$. Yet for abelian $\rho'$, $\rho'(\ell_K) = [\pm I_2]$. This makes $I_{r,1} = I_{r,0}$, and pulling back, we conclude $I^{Oh}_{r,1}|_{X_0} = I^{Oh}_{r,0}|_{X_0}$ but Equation (\ref{eq: ohtsuki's pole order formula}) suggests
\[
	{\rm pole}_{Q}(I^{Oh}_{r,1}) = 0 \text{ and } {\rm pole}_{Q}(I^{Oh}_{r,0}) \neq 0.
\]
This provides a contradiction.

As a result, $V = \hat{\tau}^{-1}(D^{irr})$ is a non-empty open subset of $X_0$. \cite[Lemma 3.15]{heusenerPorti2004PSL2Characters} then shows that $\hat{\tau}$ is injective on $V$. To summarize, we obtain a rational map $\hat{\tau}: X_0 \to D$ which is injective on a non-empty Zariski-open subset. A standard argument shows that $\hat{\tau}$ is birational. See, for example, \cite[Exercise 7.8 and paragraph right before Exercise 7.10]{harris1992algebraic}.

As a consequence, $\widetilde{X}_0$ also serves as a smooth projective model for $D$. Write the birational map $\tilde{\phi}_D: \widetilde{X}_0 \dashrightarrow D$. We now argue that $Q$ is also an ideal point for $D$: if not, $Q\in {\rm Dom}(\tilde{\phi}_D)$. Yet $I_{1,0}$ is a regular function on all of $X_{\PSL_2(\mathbb{C})}(K)$, which means it is regular on a neighborhood of $\tilde{\phi}_{D}(Q)$ and its pullback $\tilde{\phi}_D^*I_{1,0}$ is regular in a neighborhood of $Q\in \widetilde{X}_0$: this implies that ${\rm pole}_{Q}(I_{1,0}) = 0$, which contradicts Equation \eqref{eq: ohtsuki's pole order formula}. Thus, the curve $D$ and its ideal point $Q$ are the desired data on $X_{\PSL_2(\mathbb{C})}(K)$, i.e., for any $(p,q)\in \mathbb{Z}^2$, we have
\[
	{\rm pole}_Q(I_{p,q}|_D) = \varepsilon'\cdot |p - rq|.
\]

 We can then lift $D$ to the desired curve $C'$ in $X_{\SL_2(\mathbb{C})}(K)$ via Corollary \ref{cor: C to D}. The ideal point and associated valuation are liftable as well by Lemma \ref{lem: surjective morphism on smooth projective model} and Lemma \ref{lem: valuation and surjective morphism}. This completes the proof.
\end{proof}

For a 2-bridge knot $K$, let $C'\subset X_{\SL_2(\mathbb{C})}(K)$ be the curve supplied by Proposition \ref{prop: 2bridge, valuation on I_p,q}, and use the fixed peripheral system $(m_K,\ell_K)$ . We now proceed to construct the desired curve $C\subset T$ (with iterable branch) to establish Theorem \ref{thm: 2-bridge knot admits iterable branch}. The basic idea is to twist an $\SL_2(\mathbb{C})$ representation to obtain a $\mu\operatorname{rk}$-one representation, though a detailed construction shall pass to what we call the eigenvalue cover $E$ of $C'$. A conceptual diagram is the following.
\[
C' \xlongleftarrow{\pi} E \xlongrightarrow{\zeta} C\subset T.
\]

Next, we construct $E$. First, Culler--Shalen \cite[Proof of Corollary~1.4.5]{cullerShalen1983varieties} show that characters from irreducible representations form a Zariski-open subset of $X_{\SL_2(\mathbb{C})}(K)$. Hence on $C'$ they form a possibly empty open subset $U$. We argue that this is nonempty. Otherwise all characters come from reducible $\SL_2(\mathbb{C})$ representations. Such a representation can be written in upper-triangular form, and its diagonal part is an abelian representation with the same character. Yet for an abelian representation $\ell_K$ is sent to $I_2$, so $I_{r,1}=I_{r,0}$, violating the pole order formula in Proposition \ref{prop: 2bridge, valuation on I_p,q}. Thus $U\subset C'$ is a non-empty and hence dense open subset of $C'$.

Define $X\subset C'\times (\mathbb{C}^*)^2$ by the equations
\begin{equation}\label{eq: def eq for eigenvalue cover}
\begin{aligned}
M^2 - t_{m_K}M + 1 &= 0\\
L^2 - t_{\ell_K}L + 1 & = 0\\
(ML)^2 -t_{m_K\ell_K}(ML) + 1 &= 0
\end{aligned}
\end{equation}

$X$ is an affine algebraic set, since the trace functions are polynomial functions on $X_{\SL_2(\mathbb{C})}(K)$ and hence on $C'$. The projection $\pi: X \to C'$ is surjective: $m_K$ and $\ell_K$ commute in $\pi_K$ by assumption, so for a representation $\rho$ we can always take the eigenvalues of $\rho(m_K)$ and $\rho({\ell_K})$ on a common eigenline to satisfy Equation~\eqref{eq: def eq for eigenvalue cover}. Note the third equation requires the compatibility between $M$ and $L$. The defining equations also imply that $\mathbb{C}[X]$ is a finitely generated module over $\mathbb{C}[C']$ via the pullback $\pi^*$: four generators are enough: $1, M, L, ML$. As a result, the proof of Corollary~\ref{cor: C to D} applies again: there exists an irreducible Zariski-closed subset $E\subset X$ such that $\pi(E) = C'$. $\mathbb{C}[E]$ is a quotient of $\mathbb{C}[X]$, so it is also a finitely generated module over $\mathbb{C}[C']$ via pullback. Then $E$ is itself an algebraic curve. We have the following property of $E$.

\begin{lemma}\label{lem: evaluation of M and L}
	Let $K\subset S^3$ be a $2$-bridge knot. Let $r$ be a non-zero even boundary slope arising from Lemma \ref{lem: non-zero even boundary slope}, $C'\subset X_{\SL_2(\mathbb{C})}(K)$ be the algebraic curve arising from Proposition \ref{prop: 2bridge, valuation on I_p,q} and $\pi: E\to C'$ be the eigenvalue cover of $C'$ together with regular functions $M$ and $L$ as above. Then $E$ admits an ideal point $Q'$ such that 
\begin{equation}\label{eq: evaluation of M and L}
\nu_{Q'}(M) \neq 0,~r\nu_{Q'}(M) + \nu_{Q'}(L) = 0.
\end{equation}
\end{lemma}

\begin{proof}
Let $(\widetilde{E}, \tilde{\phi}_E)$ be a smooth projective model of $E$ as in Definition~\ref{defn: smooth projective model}. From Lemma~\ref{lem: surjective morphism on smooth projective model}, we obtain a compatible surjective morphism $\widetilde{\pi}: \widetilde {E} \to \widetilde{C}'$. Take $Q'\in \widetilde{\pi}^{-1}(Q)$. We argue that $Q'$ must be an ideal point of $E$. Otherwise, $\tilde{\phi}_E(Q')\in E$. The function $I_{1,0}$ is regular on $C'$, and hence $\pi^{*}I_{1,0}$ is regular on $E$. Further pulling it back, $\tilde{\phi}_E^*\pi^{*}I_{1,0}$ is a regular function on a Zariski-open neighborhood of $Q'$ so it has pole order $0$. On the other hand, Lemma~\ref{lem: valuation and surjective morphism} implies that 

\begin{equation}\label{eq: pole order on eigenvalue cover}
\operatorname{pole}_{Q'}(\pi^*I_{p,q}) = \operatorname{pole}_{Q'}(\widetilde{\pi}^*I_{p,q}) = \varepsilon_1\cdot |p - rq|,~\varepsilon_1 > 0.
\end{equation}

This provides a contradiction to the compatibility of $\widetilde{\pi}$ by taking $p=1$ and $q=0$. Hence we conclude that $Q'$ remains an ideal point of $E$. Now we further explore $\nu_{Q'}(M)$ and $\nu_{Q'}(L)$.

Taking $p=1$ and $q=0$ in Equation~\eqref{eq: pole order on eigenvalue cover}, together with $\pi^*I_{1,0}=(M-M^{-1})^2$, shows that $\nu_{Q'}(M)\neq0$. Hence
\[
 E_0=\{e\in E:M(e)^2\neq1\}
\]
is a nonempty open subset of $E$.

The third equation in \eqref{eq: def eq for eigenvalue cover} ensures
that $M$ and $L$ belong to the same peripheral eigenline on $E_0$.
To see this, fix $e\in E_0$ and a representation of the character
$\pi(e)$. Write $\alpha=M(e)$. In a basis of meridian eigenvectors,
the meridian and longitude matrices are
$\operatorname{diag}(\alpha,\alpha^{-1})$ and
$\operatorname{diag}(\beta,\beta^{-1})$, respectively.
The second equation allows $L(e)=\beta$ or $\beta^{-1}$.
In the latter case, the third equation implies
\[
 (\alpha-\alpha^{-1})(\beta-\beta^{-1})=0.
\]
Since $\alpha^2\neq1$, this forces $\beta=\beta^{-1}$, so either choice
is compatible. It follows on $E_0$, and hence on $E$, that
\[
 \pi^*I_{p,q}=(M^pL^q-M^{-p}L^{-q})^2
 \qquad ((p,q)\in\mathbb Z^2).
\]

Taking $p=r$ and $q=1$ in Equation~\eqref{eq: pole order on eigenvalue cover} shows that $r\nu_{Q'}(M)+\nu_{Q'}(L)=0$. Thus we obtain the desired relation \eqref{eq: evaluation of M and L}.

\end{proof}

\begin{proof}[Proof of Theorem \ref{thm: 2-bridge knot admits iterable branch}]
	Let $K$ be a 2-bridge knot, $r$ be a non-zero even integral boundary slope as in Lemma \ref{lem: non-zero even boundary slope}, $(m_K,\ell_K)$ be chosen peripheral data, and $C'\subset X_{\SL_2(\mathbb{C})}(K)$ be an irreducible algebraic curve with ideal point $Q$ with special pole order formula supplied by Proposition \ref{prop: 2bridge, valuation on I_p,q}. Let $\pi:E\to C'$ be the eigenvalue cover constructed as above, together with regular functions $M$ and $L$, and $E$'s smooth projective model $(\widetilde{E}, \tilde{\phi}_E)$ with ideal point $Q'\in\widetilde{E}$ and valuation data as in Lemma \ref{lem: evaluation of M and L}.
	
	We now construct the actual rank-$2$ augmentation peripheral curve $C\subset T$ as follows. Take $\sigma = 1$ if $\nu_{Q'}(M) < 0$ and $\sigma = -1$ otherwise. We define a map $\zeta: E \to T = (\mathbb{C}^*)_{(\mu,\eta)}^2$ by $\mu = M^{2\sigma}$ and $\eta = L^{2\sigma}$. Let $C\subset T$ be the Zariski closure of $\im\zeta$. Since $E$ is irreducible, so is $C$. The dimension theorem of Ferrer Santos--Rittatore \cite[Chapter~1, Theorem~4.84]{ferrersantos2005actions1} then implies that $C$ has dimension at most $1$. Dimension $0$ implies that $C$ must be a single point (see Section~\ref{subsec:iteration}). If $C$ were a point, then the map $\zeta$ would be constant. Yet (viewing $\mu$ as a regular function on $E$) $\nu_{Q'}(\mu) = 2\sigma\cdot \nu_{Q'}(M) < 0$ so $\zeta$ is not constant. Consequently, $C$ is an affine algebraic curve. Take its smooth projective model $(\widetilde{C},\tilde{\phi}_C)$. Lemma~\ref{lem: surjective morphism on smooth projective model} produces a compatible surjective morphism $\widetilde{E}\to \widetilde{C}$. Taking $P$ to be the image of $Q'$ under this morphism, Lemma~\ref{lem: valuation and surjective morphism} now implies that there exists a positive $\varepsilon$ such that

\[
-p = \nu_P(\mu) = \varepsilon\cdot2\sigma\cdot \nu_{Q'}(M) < 0.
\]

The same lemma gives
\[
 q=\nu_P(\eta)=\varepsilon\cdot2\sigma\cdot\nu_{Q'}(L),
 \qquad \frac{q}{p}= r.
\]
The pole of the regular function $\mu$ at $P$ shows that $P$ is an ideal point of $C$. Since $p>0$ and $r\notin[0,1]$, either $q>p$ or $q<0$.

As a consequence, with Lemma~\ref{lem: from algebra to analysis}, we conclude that $C$ has an iterable branch as defined in Definition~\ref{def:escaping-branch}. It remains to show that $C$ provides an actual rank-two augmentation peripheral curve, as in Definition~\ref{def:actual-kch-curve}.

First, note there is a dense open subset $U\subset C'\subset X_{\SL_2(\mathbb{C})}(K)$ consisting of characters of irreducible representations. Pulling it back via the surjective $\pi$, we obtain an open dense subset $V'=\pi^{-1}(U)\subset E$. We restrict $V'$ to the nonempty open subset where $M,L$ belong to a common peripheral eigenline, and further exclude the closed subsets $M=\pm1$. Call the resulting subset $V$. It is still open, and note that $E$ is not fully contained in $M=\pm1$ since $\nu_{Q'}(M)\neq0$ by Lemma \ref{lem: evaluation of M and L}. Hence we conclude that $V\subset E$ is an open and dense subset.

For $e\in V$, choose an irreducible representation $\rho_e$ with character $\pi(e)$, and write $\alpha=M(e)$ and $\beta=L(e)$. Let $\operatorname{ab}:\pi_K\to\mathbb Z$ be the abelianization with $\operatorname{ab}(m_K)=1$. Define the scalar twist $\rho_e^\sigma\colon\pi_K\to\GL(\mathbb C^2)$ by
\[
 \rho_e^\sigma(g)=\alpha^{\sigma\operatorname{ab}(g)}\rho_e(g).
\]
After this scalar twist, the two eigenvalues of the meridian are $1$ and $\alpha^{2\sigma}$, and hence an irreducible $\SL_2(\mathbb{C})$ representation yields an irreducible two-dimensional diagonalizable $\murk$-one representation. Since $\operatorname{ab}(\ell_K)=0$, the longitude eigenvalue on the distinguished meridian eigenline is $\beta^\sigma$. The augmentation--representation correspondence therefore realizes
\[
 (\mu(e),\eta(e))=(\alpha^{2\sigma},\beta^{2\sigma})
\]
as an actual rank-two augmentation peripheral point.

To conclude the proof, it remains to show that $\zeta(V)$ contains an open and dense subset of $C$. Now the fact that $V$ is dense in $E$ implies that $\zeta(V)$ is dense in $\zeta(E)$, and hence dense in $C = \overline{\zeta(E)}$. Thus the constructibility theorem of Ferrer Santos--Rittatore \cite[Theorem~4.91]{ferrersantos2005actions1} applies, and the proof is complete. 
\end{proof}

\begin{corollary}\label{cor:mr-two}
Let $K$ be a nontrivial knot with $\mr(K)=2$. Then, for every $N\geq0$,
\[
 \ar(\Wh^N(K))=\mr(\Wh^N(K))=b(\Wh^N(K))=2^{N+1}.
\]
\end{corollary}

\begin{proof}
The two-meridian characterization of two-bridge knots
\cite[Corollary 3.3]{BoileauZimmermann}, with the orbifold theorem
\cite{BoileauLeebPorti}, implies that $K$ is two-bridge.
Apply Theorem~\ref{thm:main-two-bridge}.
\end{proof}
% END INCLUDED SOURCE: sections/06b_two_bridge_ranks.tex
% BEGIN INCLUDED SOURCE: sections/07_ordinary_rank.tex
\section{Ordinary group rank}\label{sec:ordinary-rank}

Having proved the two-bridge result, we compare it with ordinary group
rank. Here generators are unrestricted: $\gr(K)$ is the minimum
cardinality of a generating set of $\pi_K$.

\begin{proposition}\label{prop:gr-upper}
For every knot $K$ and every $N\geq0$,
\[
 \gr(\Wh(K))\leq\gr(K)+1,\qquad
 \gr(\Wh^N(K))\leq\gr(K)+N.
\]
\end{proposition}

\begin{proof}
The relation $\ell_{\mathcal V}x_2\ell_{\mathcal V}^{-1}=x_1$ in
\eqref{eq:dutra-presentation} gives
$x_2=\ell_{\mathcal V}^{-1}x_1\ell_{\mathcal V}$. Thus $x_1$ and $\ell_{\mathcal V}$ generate the
pattern group. In the satellite group, $\ell_{\mathcal V}=\ell_K$ belongs to the
image of the companion group. It follows that $\pi_{\Wh(K)}$ is
generated by this image and $x_1$. This proves the first inequality,
and induction gives the second.
\end{proof}

For comparison, Gu--Wang--Zou \cite[Proposition 1.3]{GuWangZou}
prove that
\[
 \gr(\Wh^N(K))\geq N+1
\]
for every nontrivial knot $K$ and $N\geq1$. Each Whitehead layer
contributes a hyperbolic piece to the JSJ decomposition, and their
argument uses these pieces to bound the rank from below. Combining
this with Proposition~\ref{prop:gr-upper}, we obtain
\[
 N+1\leq\gr(\Wh^N(K))\leq N+\gr(K),\qquad N\geq1.
\]
In particular, for a fixed nontrivial $K$,
\[
 \lim_{N\to\infty}\frac{\gr(\Wh^N(K))}{N+1}=1.
\]
This asymptotic formula also appears in \cite[Corollary 2.5]{GuWangZou},
where the upper bound is expressed using tunnel number.

When $K$ is hyperbolic with tunnel number one,
\cite[Corollary 2.6]{GuWangZou} gives the exact value
$\gr(\Wh^N(K))=N+2$ for $N\geq1$. In particular, for a hyperbolic
two-bridge knot, Theorem~\ref{thm:main-two-bridge} gives the comparison
\[
 \gr(\Wh^N(K))=N+2<2^{N+1}
 =\ar(\Wh^N(K))=\mr(\Wh^N(K)),\qquad N\geq1.
\]
The exact ordinary-rank formula here uses the hyperbolic
hypothesis; we do not assert it for every two-bridge companion.

\begin{corollary}\label{cor:gr-comparison}
Suppose $K$ has a rank-$d$ augmentation $\epsilon$ with
$\epsilon(\mu)\neq1$ and peripheral point in $\mathcal S_N$. If
$2^Nd\geq\gr(K)+N$, then
\[
 \gr(\Wh^N(K))\leq\ar(\Wh^N(K)).
\]
\end{corollary}

\begin{proof}
Proposition~\ref{thm:finite-depth-iteration} gives
$\ar(\Wh^N(K))\geq2^Nd$, whereas Proposition~\ref{prop:gr-upper}
gives $\gr(\Wh^N(K))\leq\gr(K)+N$. The numerical hypothesis
compares these bounds.
\end{proof}

For a rank-two augmentation the numerical condition is
$\gr(K)\leq2^{N+1}-N$. In particular it is automatic when
$\gr(K)=2$ and $N\geq1$, provided the required peripheral point lies
in $\mathcal S_N$. The existence of such a point, for each $N$, is
supplied by an actual
rank-two augmentation peripheral curve with an iterable branch; see
the proof of Corollary~\ref{cor:escaping-rank-growth}.

\begin{proof}[Proof of Theorem~\ref{thm:main-ordinary}]
The unconditional upper bound is Proposition~\ref{prop:gr-upper}.
Under the actual-curve hypothesis,
Corollary~\ref{cor:escaping-rank-growth} gives
$\ar(\Wh^N(K))\geq2^Nd$ for every $N\geq0$.
Consequently,
\[
 \gr(\Wh^N(K))\leq\gr(K)+N\leq2^Nd\leq\ar(\Wh^N(K))
\]
whenever the middle inequality holds.

Finally, a nontrivial two-bridge knot has a two-generator nonabelian knot
group, so $\gr(K)=2$. Theorem~\ref{thm:main-two-bridge} and
$N+2\leq2^{N+1}$ yield the final assertion.
\end{proof}

\begin{remark}
For a fixed knot with the actual-curve hypothesis and $d>0$, the numerical
condition in Theorem~\ref{thm:main-ordinary} holds for all sufficiently
large $N$. This conclusion depends on the curve hypothesis; the linear
upper bound alone does not imply $\gr(\Wh^N(K))\leq\ar(\Wh^N(K))$.
\end{remark}
% END INCLUDED SOURCE: sections/07_ordinary_rank.tex
% BEGIN INCLUDED SOURCE: sections/08_first_doubles.tex

\section{The first Whitehead double of the square knot}\label{sec:first-doubles}

In this section, we prove Theorem~\ref{thm:main-square}. Let
\[
S:=T(2,3)\mathbin{\#}\overline{T(2,3)}
\]
be the square knot. It has 
\[
\ar(S)=\mr(S)=b(S)=3
\]
where $\ar(S)=2+2-1=3$ is given by the formula in \cite[Theorem~1.9]{Cornwell}, and $b(S) = 2+2-1 =3$ is by \cite{SchultensBridge}. 

For its untwisted Whitehead double $\Wh(S)$, Dutra's bridge formula gives $b(\Wh(S))=2b(S)=6$ \cite[Section~4]{Dutra}. In Section \ref{section:ar-S}, we show that $\ar(\Wh(S))\le5$, and in Section \ref{sec:square-meridional}, we then use a finite quotient and group cohomology to show
that its $\mr(\Wh(S))=6$.

\subsection{Augmentation rank of \texorpdfstring{$\Wh(S)$}{Wh(S)}}
\label{section:ar-S}

The knot group of the square knot is
\[
 \pi_S=\langle a,b,c\mid aba=bab,\ aca=cac\rangle ,
\]
where the peripheral subgroup is generated by $m = a, \ell =(ab)^3(ac)^{-3}$. See e.g. \cite[Section~2.2]{TuffleySquare} for such a presentation. By Proposition \ref{prop:van-kampen},

\begin{equation}\label{eq:square-satellite-presentation}
 \pi_{\Wh(S)}=\left\langle x_1,x_2,a,b,c,\ell\ \middle|\
 \begin{aligned}
 a&=x_2x_1^{-1}x_2^{-1}x_1, & \ell x_2\ell^{-1}&=x_1, & [a,\ell]&=1,\\
 &aba=bab, & aca=&cac, & \ell=(ab)&^3(ac)^{-3}
 \end{aligned}\right\rangle .
\end{equation}

\begin{proposition}
\label{thm:square-knot-counterexample}
Over $\C$, $\ar(\Wh(S))\le5$.

\end{proposition}

\begin{proof}
We know $\ar(\Wh(S))\leq b(\Wh(S)) =6$. We shall prove by contradiction and assume $\ar(\Wh(S))=6$. By definition, there exists a 6-dimensional irreducible diagonalizable $\murk$-1 representation
\[
 \varphi \colon\pi_{\Wh(S)}\longrightarrow \GL(V).
\]
Let $X_i=\varphi(x_i), A = \varphi(a)$ etc. Let $W_i := \im(X_i-I), W_a := \im(A-I)$ etc. Since $(\varphi, V)$ is $\murk$-1, each of $W_i$ has dimension $1$. Since $a,b,c$ are conjugate, their microlocal ranks are equal. Moreover, in the doubling construction, each of them is the composite of two meridians of $\Wh(S)$, therefore
\[
s:=\dim W_a = \dim W_b = \dim W_c \leq 2.
\]
Define $W_{12} = W_1+W_2$, and define $W_{abc} = W_a+W_b+W_c$.

We consider the subspace $W:= W_{12} + W_{abc}\subset V$. Each of $X_i$ preserves every subspace containing $W_{12}$, hence preserving $W$. The matrices $A,B,C$ preserve $W_{abc}$ and acts as the identity on $V/W_{abc}$, hence they also preserve $W$. Since $\ell$ is a word of $a,b,c$, its image preserves $W$ as well. Altogether, we see $\pi_{\Wh(S)}$ preserves $W$. By irreducibility, we have $\dim V = \dim W$. Therefore,
\[
 6=\dim V = \dim W \le\dim W_{12}+ \dim W_{abc}\le2+3s.
\]
Since $s\le2$, it is only possible with $s=2$.  Now
$W_a\subset W_{12}$ and both spaces have dimension two, so
$W_a=W_{12}$. Consequently $W$ itself is invariant under
$X_1,X_2$ as well as the companion group.  Irreducibility and dimensions
give the direct decomposition
\[
 V=W_a\oplus W_b\oplus W_c.
\]

By choosing a suitable basis, we have
\[
A=
\begin{pmatrix}
N&A_{12}&A_{13}\\
0&I&0\\
0&0&I
\end{pmatrix},\qquad
B=
\begin{pmatrix}
I&0&0\\
B_{21}&N&C_{21}\\
0&0&I
\end{pmatrix},\qquad
C=
\begin{pmatrix}
I&0&0\\
0&I&0\\
C_{31}&C_{32}&N
\end{pmatrix},
\]
where \(N\in \mathrm{GL}_2(\mathbb C)\). Direct block multiplication in \(ABA=BAB\) gives
\[
(N+A_{12}B_{21})[N-I,A_{12},A_{13}]=0,
\qquad
(B_{21}A_{12}+N)[B_{21},N-I,C_{21}]=0.
\]
The matrices $[N-I,A_{12},A_{13}], [B_{21},N-I,C_{21}]$ have rank two because $A,B$ have  $\murk$=2. Therefore
\[
A_{12}B_{21}=B_{21}A_{12}=-N.
\]
Similarly, the relation \(ACA=CAC\) gives
\[
A_{13}C_{31}=C_{31}A_{13}=-N.
\]

With these identities, it is a routine calculation to find 
\[
L|_{W_a} = (AB)^3(AC)^{-3}|_{W_a}= -N^3 (-N^3)^{-1} = I.
\]

Finally, since $A=X_2X_1^{-1}X_2^{-1}X_1$, we have $W_a\subset W_1 + W_2$. By dimension count, we must have $W_a = W_1\oplus W_2$. However, the relation \(LX_2L^{-1}=X_1\) yields
$W_2=L^{-1}W_1$. Since $L|_{W_a}=I$, this implies $W_1 = W_2$, a contradiction.

This concludes the proof by contradiction, and we prove the assertion.
\end{proof}

\subsection{Meridional rank of \texorpdfstring{$\Wh(S)$}{Wh(S)}}
\label{sec:square-meridional}

We use group cohomology to give an estimate for $\mr(\Wh(S))$. We recall the basic notions; see
\cite{Brown1994,HiltonStammbach1997} for general references. Let $G$ be a group and $V$ be a left $G$-module. Functorially, $H^n(G,V) := R^n\operatorname{Hom}_{\mathbb{Z}[G]}(\mathbb{Z},V)$. Alternatively, the bar resolution gives a concrete cochain complex
\[
C^n(G, V) = \operatorname{Map}(G^n, V),
\]
where the maps are set theoretic, and the differential is given by 
\begin{align*}
	(df)(g_1, \dotsb, g_{n+1}) = &g_1\cdot f(g_2, \dotsb, g_{n+1}) + \sum_{i=1}^{n}(-1)^i f(g_1,\dotsb, g_ig_{i+1},\dotsb, g_{n+1}) \\
	& + (-1)^{n+1}f(g_1, \dotsb, g_{n}).
\end{align*}

We note two special cases. For $n=0$,
\[
H^0(G,V) = V^G :=\{v\in V:\rho(g)v=v\text{ for every }g\in G\}.
\]
For $n=1$, the 1-cocycle maps are given by crossed homomorphisms:
\[
Z^1(G;V)=\{f:G\to V:f(gh)=f(g)+\rho(g)f(h)\text{ for all }g,h\in G\},
\]
because $(df)(g,h) = gf(h) - f(gh) +f(g) =0$. The 1-coboundary maps are given by principal crossed homomorphisms.
\[
B^1(G,V)=\{f_v:v\in V,\ f_v(g)=\rho(g)v-v\},
\]
because $(dv)(g) = gv - v$.

\begin{lemma}\label{lem:cohomological-meridional-bound}
Let $G$ be generated by $r\geq1$ conjugates of an element $m$, and let $\rho:G\to\GL(V)$ be a representation on a finite-dimensional vector space over a field $\mathbb F$. Then
\begin{equation}\label{eq:cohomological-meridional-bound}
 \dim Z^1(G,V)-\dim V \leq (r-1)\murk(\rho)
\end{equation}
\end{lemma}

\begin{proof}

For any $1$-cocycle $f:G\to V$, and any $g\in G$, apply $f$ to the identity $e_G = g g^{-1}$, we get $0= f(g) + \rho(g)f(g^{-1})$. Therefore $f(g^{-1}) = -\rho(g^{-1})f(g)$.

Let $m_1=m$ and $m_i=g_i^{-1} m g_i$ for $2\leq i\leq r$. We have
\begin{align}
f(m_i)= f (g_i^{-1} m g_i) 
 &= f(g_i^{-1}) + \rho(g_i^{-1}) f(m g_i) \\
 &= -\rho(g_i^{-1})f(g_i) + \rho(g_i^{-1}) (f(m) + \rho(m)f(g_i))\\
 &=\rho(g_i^{-1}) f(m) + \rho(g_i^{-1})(\rho(m)-I)f(g_i)\\
 &=\rho(g_i^{-1}) f(m) + (\rho(m_i)-I)\rho(g_i^{-1})f(g_i).
\end{align}

A cocycle is determined by its values on the generating set. Once $f(m)$ is fixed, the possible value of $f(m_i)$ belongs to an affine space with direction $\im((\rho(m_i)-I))$, the dimension of which equals to the microlocal rank. Therefore, the map
\begin{align*}
	Z^1(G,V)&\longrightarrow V\oplus\bigoplus_{i=2}^r\im(I-\rho(m_i)),\\
	f&\longmapsto\bigl(f(m),(f(m_i)-\rho(g_i^{-1}) f(m))_{i=2}^r\bigr),
\end{align*}
is injective. Thus $\dim Z^1(G;V)\leq\dim V+(r-1)\murk(\rho)$, the inequality follows.

On the other hand, there is a short exact sequence $0\rightarrow H^0(G;V) \rightarrow V\rightarrow B^1(G;V)\rightarrow 0$. This implies $\dim V = \dim B^1(G;V) + \dim H^0(G;V)$. 
\end{proof}

\begin{remark}
	The left hand side of the inequality is in fact the truncated negative Euler number $\dim H^1(G,V)-\dim H^0(G,V)$. Indeed, there is a short exact sequence $0\rightarrow H^0(G,V) \rightarrow V\rightarrow B^1(G,V)\rightarrow 0$. This implies $\dim V = \dim B^1(G,V) + \dim H^0(G,V)$. Therefore,
	\[
	\dim Z^1(G,V)-\dim V = \dim H^1(G,V)-\dim H^0(G,V).
	\]
\end{remark}

Let $S_{14}$ be the permutation group of 14 elements. We define a group homomorphism $\phi:\pi_{\Wh(S)}\to S_{14}$ by 
\begin{equation}\label{eq:square-s14-generators}
\begin{aligned}
\phi(x_1)&=(1\,2)(3\,4)(5\,6),\\
\phi(x_2)&=(1\,7)(3\,8)(5\,9),\\
\phi(a)&=(1\,7\,2)(3\,8\,4)(5\,9\,6),\\
\phi(b)&=(1\,13\,7)(3\,5\,10)(4\,11\,9),\\
\phi(c)&=(1\,3\,10)(2\,12\,8)(5\,14\,9),\\
\phi(\ell)&=(1\,6\,8)(2\,9\,3)(4\,7\,5)(10\,12\,11)
\end{aligned}
\end{equation}
It is straightforward to check these assignments satisfy the relators in $\pi_{\Wh(S)}$, hence $\phi$ is well defined.

\begin{proposition}\label{prop:square-s14}
$\phi$ is surjective.
\end{proposition}

\begin{proof}
Observe that
\[
\phi(x_1bc)=(1\,6\,5\,14\,3\,4\,11\,9\,10\,13\,7\,2\,12\,8)
\]
is a cyclic permutation of order $14$, and 
\[
\phi((x_2b)^7)=(1\,13)
\]
is a simple transposition. The distance between 1 and 13 in the cyclic permutation is 9, which is coprime to 14. Therefore these two elements generate $S_{14}$, proving surjectivity.
\end{proof}

Now let $V = \mathbb F_2^{14}$. Define the representation $\rho_\phi$ by the composition
\[
\rho_\phi: \pi_{\Wh(S)}\xrightarrow{\phi} S_{14} \rightarrow GL(V),
\]
where the second map is the permutation representation, defined by permuting the basis vectors.

\begin{proposition}\label{prop:square-cohomology}
$(\rho_\phi, V)$ has microlocal rank $3$ and 
\begin{equation}\label{eq:square-cohomology-dimensions}
 \dim Z^1(\pi_{\Wh(S)}, V)=27.
\end{equation}
\end{proposition}

\begin{proof}
We compute by brute force. Let $e_i = (0,\ldots,0,1,0,\ldots,0)$ where the $i$-th coordinate is 1. Then the standard basis for $V$ is $\{e_1, \ldots, e_{14}\}$. Since $\phi(x_1)=(1\,2)(3\,4)(5\,6)$, three simple transpositions. We can compute 
\[
\im (\rho_\phi(x_1)- I) = \operatorname{Span}\{e_1+e_2, e_3+e_4, e_5+e_6\}.
\]
Hence, the microlocal rank is 3.

For the cocycles, recall that $\pi_{\Wh(S)}$ is generated by 6 generators, relabeled as $s_i$ for $i=1,\ldots,6$, and 6 relators, denoted by $R_i$. To name a function $f$, we need to specify $f(s_i)$ for all $i=1,\ldots,6$. Each $f(s_i)$ is a linear combination of 14 basis vectors. Therefore, $f$ is specified by a vector of $84 = 6 \times 14$ binary coordinates. We can choose $\{(s_i,e_j)_{1\leq i\leq 6, 1\leq j\leq 14}\}$ to be the basis of this vector space $U$.

To be well-defined, $f$ must satisfy the cocycle condition for each relator $R_i$. For example, take $R_2$ to be $e_{\pi_{\Wh(S)}} =\ell x_2 \ell^{-1}x_1^{-1}$. Apply $f$ to both sides, we get
Applying the cocycle \(f\) to
\(e_{\pi_{\operatorname{Wh}(S)}}=\ell x_2\ell^{-1}x_1^{-1}\), we obtain
\begin{align*}
0
&=f(\ell)
+\rho_\phi(\ell)f(x_2)
+\rho_\phi(\ell x_2)f(\ell^{-1})
+\rho_\phi(\ell x_2\ell^{-1})f(x_1^{-1})\\
&=f(\ell)+\rho_\phi(\ell)f(x_2)
-\rho_\phi(\ell x_2\ell^{-1})f(\ell)
-\rho_\phi(\ell x_2\ell^{-1}x_1^{-1})f(x_1)\\
&=\bigl(I-\rho_\phi(x_1)\bigr)f(\ell)
+\rho_\phi(\ell)f(x_2)-f(x_1).
\end{align*}
This is a system of $\dim V =14$ linear equations. Solutions to this system is the kernel of a $14 \times 84$ matrix. Since there are 6 relators, we obtain a total of $14 \times 6=84$ equations, which form a matrix of size $84 \times 84$.

For clarity, we describe how the matrix is assembled. Let
\[
(s_1,\ldots,s_6)=(x_1,x_2,a,b,c,\ell),
\qquad M_i=\rho_\phi(s_i),
\]
where $M_i e_j=e_{\phi(s_i)(j)}$. We write the six relators as
\[
\begin{aligned}
R_1&=a^{-1}x_2x_1^{-1}x_2^{-1}x_1,\\
R_2&=\ell x_2\ell^{-1}x_1^{-1},\\
R_3&=a\ell a^{-1}\ell^{-1},\\
R_4&=aba b^{-1}a^{-1}b^{-1},\\
R_5&=aca c^{-1}a^{-1}c^{-1},\\
R_6&=\ell^{-1}(ab)^3(c^{-1}a^{-1})^3.
\end{aligned}
\]
In the following pseudocode, $D[r,i]$ denotes the
$14\times14$ block in block row $r$ and block column $i$.
All additions are over $\mathbb F_2$, and powers in the
relators are expanded before reading the letters.

\begin{verbatim}
D := zero matrix of size 84 x 84
for r = 1,...,6:
    P := identity matrix of size 14
    for each letter s_i^epsilon in R_r, from left to right:
        if epsilon = -1:
            P := P M_i^(-1)
        D[r,i] := D[r,i] + P
        if epsilon = 1:
            P := P M_i
return D
\end{verbatim}

At the start of each letter, $P$ is the image of the
preceding prefix. The order of the updates follows from
the cocycle identity and
$f(s_i^{-1})=-M_i^{-1}f(s_i)$; the minus sign disappears
over $\mathbb F_2$. Thus the $r$-th block row records
$f(R_r)=0$. This gives the matrix displayed below.

\begin{center}
\begingroup
\fontsize{5}{6}\selectfont
\setlength{\arraycolsep}{0.5pt}
\renewcommand{\arraystretch}{0.9}

\begin{adjustbox}{
  max totalsize={\linewidth}{0.9\textheight},
  center
}
\(
\left(
\begin{array}{*{84}{c}}
\Z{15}1 & \Z{4}1 & \Z{13}1 & \Z{48}0 \\
1 & \Z{5}1 & \Z{21}1 & \Z{54}0 \\
\Z{17}1 & \Z{3}1 & \Z{13}1 & \Z{47}0 \\
\Z{2}1 & \Z{4}1 & \Z{22}1 & \Z{52}0 \\
\Z{19}1 & \Z{2}1 & \Z{13}1 & \Z{46}0 \\
\Z{4}1 & \Z{3}1 & \Z{23}1 & \Z{50}0 \\
1 & \Z{5}1 & \Z{8}1 & \Z{4}1 & \Z{8}1 & \Z{53}0 \\
\Z{2}1 & \Z{4}1 & \Z{9}1 & \Z{3}1 & \Z{9}1 & \Z{51}0 \\
\Z{4}1 & \Z{3}1 & \Z{10}1 & \Z{2}1 & \Z{10}1 & \Z{49}0 \\
\Z{37}1 & \Z{45}0 \\
\Z{38}1 & \Z{44}0 \\
\Z{39}1 & \Z{43}0 \\
\Z{40}1 & \Z{42}0 \\
\Z{41}1 & \Z{41}0 \\
1 & \Z{20}1 & \Z{48}1 & 1 & \Z{11}0 \\
\Z{1}1 & \Z{14}1 & \Z{53}1 & 1 & \Z{11}0 \\
\Z{2}1 & \Z{19}1 & \Z{49}1 & 1 & \Z{9}0 \\
\Z{3}1 & \Z{14}1 & \Z{53}1 & 1 & \Z{9}0 \\
\Z{4}1 & \Z{15}1 & \Z{53}1 & 1 & \Z{7}0 \\
\Z{5}1 & \Z{8}1 & \Z{59}1 & 1 & \Z{7}0 \\
\Z{6}1 & \Z{10}1 & \Z{65}0 \\
\Z{7}1 & \Z{11}1 & \Z{63}0 \\
\Z{8}1 & \Z{6}1 & \Z{67}0 \\
\Z{9}1 & \Z{14}1 & \Z{58}0 \\
\Z{10}1 & \Z{14}1 & \Z{57}0 \\
\Z{11}1 & \Z{11}1 & \Z{59}0 \\
\Z{12}1 & \Z{13}1 & \Z{56}0 \\
\Z{13}1 & \Z{13}1 & \Z{55}0 \\
\Z{28}1 & \Z{6}1 & \Z{34}1 & 1 & \Z{11}0 \\
\Z{29}1 & 1 & \Z{40}1 & \Z{4}1 & \Z{6}0 \\
\Z{30}1 & \Z{5}1 & \Z{35}1 & 1 & \Z{9}0 \\
\Z{31}1 & 1 & \Z{40}1 & \Z{3}1 & \Z{5}0 \\
\Z{32}1 & \Z{1}1 & \Z{39}1 & 1 & \Z{7}0 \\
\Z{28}1 & \Z{4}1 & \Z{41}1 & \Z{2}1 & \Z{4}0 \\
\Z{31}1 & \Z{2}1 & \Z{35}1 & \Z{5}1 & \Z{6}0 \\
\Z{33}1 & \Z{1}1 & \Z{36}1 & \Z{4}1 & \Z{5}0 \\
\Z{29}1 & \Z{6}1 & \Z{37}1 & \Z{3}1 & \Z{4}0 \\
\Z{37}1 & 1 & \Z{44}0 \\
\Z{38}1 & 1 & \Z{43}0 \\
\Z{37}1 & \Z{1}1 & \Z{43}0 \\
\Z{83}0 \\
\Z{83}0 \\
\Z{28}1 & 1 & \Z{4}1 & \Z{8}1 & \Z{39}0 \\
\Z{40}1 & \Z{2}1 & \Z{39}0 \\
\Z{30}1 & \Z{5}1 & 1 & \Z{6}1 & 1 & \Z{5}1 & \Z{31}0 \\
\Z{31}1 & \Z{3}1 & 1 & \Z{8}1 & 1 & \Z{2}1 & \Z{33}0 \\
\Z{30}1 & \Z{1}1 & 1 & \Z{11}1 & 1 & 1 & \Z{35}0 \\
\Z{38}1 & \Z{8}1 & \Z{35}0 \\
\Z{40}1 & \Z{1}1 & \Z{5}1 & \Z{5}1 & \Z{28}0 \\
\Z{37}1 & \Z{11}1 & \Z{33}0 \\
\Z{30}1 & \Z{5}1 & \Z{1}1 & \Z{7}1 & \Z{3}1 & \Z{1}1 & \Z{30}0 \\
\Z{37}1 & \Z{9}1 & \Z{35}0 \\
\Z{38}1 & \Z{10}1 & \Z{33}0 \\
\Z{39}1 & \Z{13}1 & \Z{29}0 \\
\Z{40}1 & \Z{2}1 & \Z{39}0 \\
\Z{41}1 & \Z{13}1 & \Z{27}0 \\
\Z{28}1 & \Z{6}1 & \Z{1}1 & \Z{18}1 & 1 & \Z{7}1 & \Z{17}0 \\
\Z{29}1 & \Z{4}1 & 1 & \Z{21}1 & 1 & \Z{3}1 & \Z{20}0 \\
\Z{28}1 & \Z{1}1 & 1 & \Z{25}1 & 1 & 1 & \Z{23}0 \\
\Z{39}1 & \Z{19}1 & \Z{23}0 \\
\Z{32}1 & 1 & \Z{2}1 & \Z{24}1 & \Z{21}0 \\
\Z{41}1 & \Z{19}1 & \Z{21}0 \\
\Z{37}1 & \Z{24}1 & \Z{20}0 \\
\Z{28}1 & \Z{6}1 & \Z{3}1 & \Z{18}1 & \Z{4}1 & \Z{3}1 & \Z{15}0 \\
\Z{41}1 & \Z{18}1 & \Z{3}1 & \Z{4}1 & \Z{13}0 \\
\Z{37}1 & \Z{21}1 & \Z{23}0 \\
\Z{38}1 & \Z{27}1 & \Z{16}0 \\
\Z{39}1 & \Z{22}1 & \Z{20}0 \\
\Z{40}1 & \Z{27}1 & \Z{14}0 \\
\Z{41}1 & \Z{19}1 & \Z{21}0 \\
\Z{31}1 & \Z{1}1 & \Z{1}1 & \Z{2}1 & \Z{10}1 & 1 & \Z{1}1 & \Z{5}1 & \Z{16}1 & \Z{7}0 \\
\Z{29}1 & 1 & \Z{3}1 & \Z{2}1 & \Z{7}1 & \Z{10}1 & \Z{5}1 & \Z{2}1 & \Z{12}1 & \Z{4}0 \\
\Z{28}1 & 1 & 1 & 1 & \Z{7}1 & 1 & \Z{2}1 & \Z{4}1 & \Z{5}1 & \Z{4}1 & \Z{3}1 & \Z{3}1 & \Z{3}1 & \Z{11}0 \\
\Z{29}1 & \Z{1}1 & \Z{2}1 & \Z{4}1 & \Z{2}1 & \Z{19}1 & 1 & \Z{3}1 & \Z{8}1 & \Z{6}0 \\
\Z{31}1 & 1 & 1 & \Z{1}1 & \Z{1}1 & \Z{3}1 & \Z{2}1 & \Z{4}1 & \Z{1}1 & \Z{9}1 & \Z{2}1 & \Z{4}1 & \Z{3}1 & \Z{9}0 \\
\Z{33}1 & \Z{1}1 & \Z{1}1 & \Z{3}1 & \Z{2}1 & \Z{2}1 & \Z{3}1 & \Z{9}1 & \Z{2}1 & \Z{4}1 & \Z{7}1 & \Z{5}0 \\
\Z{30}1 & \Z{1}1 & 1 & 1 & \Z{2}1 & 1 & \Z{8}1 & \Z{2}1 & \Z{1}1 & \Z{3}1 & \Z{2}1 & \Z{5}1 & \Z{8}1 & \Z{8}0 \\
\Z{29}1 & \Z{10}1 & \Z{2}1 & \Z{4}1 & \Z{5}1 & \Z{2}1 & \Z{12}1 & \Z{12}0 \\
\Z{46}1 & \Z{13}1 & \Z{11}1 & \Z{10}0 \\
\Z{30}1 & 1 & \Z{5}1 & \Z{1}1 & \Z{13}1 & \Z{5}1 & \Z{3}1 & \Z{1}1 & \Z{15}1 & \Z{1}0 \\
\Z{32}1 & 1 & \Z{3}1 & 1 & \Z{8}1 & \Z{2}1 & 1 & \Z{14}1 & \Z{12}1 & \Z{3}0 \\
\Z{29}1 & \Z{1}1 & \Z{2}1 & 1 & \Z{2}1 & 1 & \Z{4}1 & \Z{4}1 & \Z{2}1 & \Z{3}1 & \Z{5}1 & \Z{4}1 & \Z{12}1 & \Z{2}0 \\
\Z{28}1 & 1 & \Z{13}1 & \Z{4}1 & \Z{5}1 & \Z{13}1 & \Z{13}1 & 0 \\
\Z{32}1 & 1 & \Z{21}1 & \Z{5}1 & \Z{2}1 & \Z{4}1 & \Z{13}1
\end{array}
\right)
\)
\end{adjustbox}

\endgroup
\end{center}

The set of cocycles  $Z^1(\pi_{Wh(S)},V)$ is the kernel of this matrix. Using a computer program, we find that
\[
\dim Z^1(\pi_{Wh(S)},V)=27.
\]
\end{proof}

\begin{proof}[Proof of Theorem~\ref{thm:main-square}]
We know that $b(\Wh(S))=6$, and that $\ar(\Wh(S))\leq 5$ by Proposition~\ref{thm:square-knot-counterexample}. If $\pi_{Wh(S)}$ is generated by $r$ meridians, Lemma~\ref{lem:cohomological-meridional-bound} and
Proposition~\ref{prop:square-cohomology} imply
\[
 27-14 = 13\leq 3(r-1).
\]
Hence $r > 5$. Together with the upper bound $b = 6$, this implies $\mr(\Wh(S)) = 6$.
\end{proof}

\subsection*{Acknowledgements}

The authors would like to thank Yu Zhao for some helpful conversations.

The authors acknowledge the use of OpenAI's ChatGPT during
the development of this work. This project began in October 2025, motivated in part by an interest in exploring how AI tools could support mathematical research. Initially, GPT assisted with writing SageMath code, implementing computations in SnapPy, and conducting computational experiments. It was later used to explore proof strategies, develop arguments, prepare an initial draft, and assist in checking the code, computations, and proposed proofs.

Throughout this process, the authors supplied the relevant literature, many of the initial mathematical ideas, detailed prompts, and directed the successive stages of experimentation, proof development, and revision. The initial AI-assisted draft was subsequently substantially rewritten by the authors, who reorganized the manuscript and revised the exposition and presentation of the proofs. Further checks and revisions were carried out through mathematical discussions among the authors and additional interactions with these tools. The AI-assisted checks were part of this review process, not an independent certification of the proofs. The authors take full responsibility for the mathematical content, computations, references, and final presentation. 

HG is supported by the National Key R\&D Program of China 2023YFA1010500 and NSFC Grant No.~12625107.
SG is supported by NSFC Grant No.~12201102.

% END INCLUDED SOURCE: sections/08_first_doubles.tex
\printbibliography
\end{document}